\documentclass[12pt]{article}

\usepackage{amsmath,amssymb,amsthm,graphicx,hyperref,color,CJK}
\usepackage[utf8]{inputenc}
\usepackage[english]{babel}
\usepackage{mathrsfs}
\usepackage{subcaption}
\usepackage{sectsty}
\usepackage[export]{adjustbox}
\usepackage[edges]{forest}
\usepackage{xcolor}
\usepackage{titlesec}
\usepackage{comment}

\usepackage[numbers,sort&compress]{natbib}

\sectionfont{\centering}

\titleformat{\section}{\large\bfseries}{\thesection}{1em}{}

\newtheorem{theorem}{Theorem}[section]
\newtheorem{corollary}{Corollary}[theorem]
\newtheorem{lemma}[theorem]{Lemma}
\newcommand{\Z}{\mathbb{Z}}
\newcommand{\R}{\mathbb{R}}

\newcommand{\cS}{\mathcal{S}}
\newcommand {\calS}{\mathcal{S}}

\begin{document}
\vspace{-2mm}
\title{Orthogonal Splines Associated to Root Systems and Applications to Function Approximation From Generalized Derivative Samples}
\author{Dionel Jaime}
\date{}
\maketitle

\vspace{-6mm}
\begin{abstract}

Let $\mathcal{S}(X)$ be the (algebraic) principal shift invariant space generated by a box spline $B_X$ with $d$-dimensional direction set $X$ and associated differential operator $D_X$. For the special case when the directions in $X$ each have even multiplicity i.e $X = [X_0|X_0]$ and are appropriately selected from a root system, we may endow $\cS(X)$ with a special inner product $\langle \cdot , \cdot \rangle_X$ which is isometrically isomorphic to a space of algebraic characters studied by Heckman-Opdam. Using this algebraic theory, we develop novel biorthogonal systems for $\cS(X)$ which satisfy that, within any fixed biorthogonal system, the diameter of the supports of the basis functions is not uniformly bounded; however, the image of any biorthogonal system under $D_X$ (which is necessarily finite and discrete) is uniformly finite. In other words, the basis elements, while having non-sparse box-spline representation, are uniformly sparse in their truncated-power representation. \\
 
We take a biorthogonal system and construct approximation maps into $\mathcal{S}(X)$. We construct the analysis-synthesis operator from smooth compactly supported continuous functions into $\mathcal{S}(X)$ using the Sobolev-like inner product $\langle D_{X_0}\cdot, D_{X_0}\cdot \rangle_{L^2(\R^d)}$. We prove that this map is a local discrete quasi-interpolation operator and, when $X$ comes from an irreducible root system, this quasi-interpolant can be modified to yield optimal approximation order. In the case where $X$ comes from the $A_n$ root system, we extend the domain of $\langle \cdot, \cdot\rangle_X$ to produce a projector into $\mathcal{S}(X)$ which has optimal approximation order. We use integration by parts to represent these optimal maps in terms of generalized derivative samples induced by one of our bases. 
\end{abstract}
\vspace{-4mm}

\section{Introduction}

\subsection{The Type of Results Achieved In This Paper}
In order to put the content of this paper in context, we first consider a setup that is simpler than the one explored in this paper but gives an idea of the results achieved. \\

Suppose $D$ is a homogeneous differential operator acting on a Sobolev space $\mathcal{H}$. Suppose $\mathcal{S}$ is a subspace of $\mathcal{H}$ and for any $f,g \in \mathcal{S}$ the following rule defines an inner product on $\mathcal{S}$: 

\begin{equation}\label{obvious_IP}
   \langle f,g \rangle_{\mathcal{S}} := \langle Df , Dg\rangle_{L^2(\R^d)}.
\end{equation}

Suppose that $\mathcal{B}$ is an orthonormal basis for $\mathcal{S}$ under $\langle \cdot, \cdot \rangle_{\mathcal{S}}$ consisting of compactly supported functions which satisfy the following properties: 

\begin{itemize}
   \item If $f \in \mathcal{H}$ is compactly supported, $\displaystyle{ \langle Df , DE\rangle_{L^2(\R^d)} = 0 }$ for all but finitely many $E \in \mathcal{B}$.
   \item The diameters of the sets in $\{\text{supp}(D^2E) \ | \ E \in \mathcal{B} \}$ are uniformly bounded.
\end{itemize}


As a result of the first property, the following map

\begin{equation}\label{approximation_map}
   P: f \mapsto \sum_{E \in \mathcal{B}} \langle Df , DE \rangle_{L^2(\R^d)} E
\end{equation}

is an orthogonal projector from smooth compactly supported function onto $\mathcal{S}$. If $P(f)$ approximates $f$ in some sense, then integration by parts on the coefficient of each $E \in \mathcal{B}$ above implies that, $P(f)$ approximates $f$ from the collection of \textit{generalized derivative samples} induced by $\mathcal{B}$ (see e.g \cite{UnserAldroubi1994}, \cite{GarciaPerezVillalon2008} for generalized sampling terminology): 




$$\{\langle D^{2-i}f, D^iE \rangle_{L^2(\R^d)}  \}_{E \in \mathcal{B}} $$

where $0 \leq i \leq 2$. The second property satisfied by $\mathcal{B}$ implies that, while $\mathcal{B}$ is not necessarily a local basis, each basis element is the convolution of the local collection $\{D^2E| \ E \in \mathcal{B}\}$ with a single suitable Green's function i.e for all $E \in \mathcal{B}$ there exists a function $G$ such that 
$$G \ast D^2E = E. $$

In this paper, $\mathcal{S}$ is a box spline space associated to a root system. The collection of splines used to construct $P$ will be an orthonormal basis, or a biorthogonal system more generally, for $\mathcal{S}$ with respect to an inner product similar to, but different from, $\langle \cdot, \cdot \rangle_S$. To describe our set-up exactly, we review box spline basics and introduce this inner product.



\subsection{The Box Spline Space $\mathcal{S}(X)$ and The Inner Product $\langle \cdot , \cdot \rangle_X$} 

Suppose $\Gamma \subset \R^d $ is a Euclidean lattice. We fix the lattice $\Gamma$. All objects defined below may be dependent on this fixed $\Gamma$. Whenever it makes sense, we define the \textit{semi-discrete convolution} between a function $B: \R^d \to \R$ and a sequence $a: \Gamma \to \R$ as:
\begin{equation}
    B \ast' a := \sum_{\gamma \in \Gamma} a(\gamma)B(\cdot - \gamma).
\end{equation}
In lieu with this, we define the space   
\begin{equation}
    S(\Gamma) := \{a: \Gamma \to \R | \  \text{$a$ is finitely supported} \}.
\end{equation}

On $S(\Gamma)$, $f \ast' g$ denotes standard convolution and the inner product on $S(\Gamma)$ is:

\begin{equation}
    \langle f,g \rangle_{\ell^2(\Gamma)} : = \sum_{\gamma \in \Gamma} f(\gamma)g(\gamma). 
\end{equation}

The semi-discrete convolution appears in the construction of \textit{shift invariant (SI) spaces}, i.e. functions spaces on $\R^d$ that are invariant under translations by $\Gamma$. 
We refer to the image of $S(\Gamma)$ under the map $B \ast'$ as the \textit{principal shift invariant (PSI) space generated by $B$.} There is a rich history of using SI spaces to approximate subspaces of $L^2(\R^d)$ (see, e.g. \cite{Approx1} and \cite{Approx2}). PSI spaces are central to this theory, as many approximation properties of a general SI space are determined by approximation properties of the $L^2$ closure of its PSI subspaces. \\
 
The subject of this article is an important class of PSI spaces: The \textit{box spline spaces}. To define a box spline, we start with a finite multi-set of $n$ nonzero vectors $X \subset \R^d$, referred to as the \textit{direction set} or, once identifying $X$ with a matrix whose columns comprise the given multi-set, as a \textit{direction matrix}. The box spline $B_X$ is the distribution on $C(\R^d)$   

\begin{equation}\label{Box Spline Distribution Definition}
    B_X: F \mapsto \int_{[-\frac{1}{2}, \frac{1}{2}  ]^n} F(Xt) \ dt.
\end{equation}





$B_X$ is supported on the \textit{centered zonotope} associated to $X$: 

\begin{equation}
   Z(X):= \bigg\{ \sum_{v \in X} c_v v \bigg| \ -\frac{1}{2} \leq c_v \leq \frac{1}{2} \ \text{for all} \ v \in X  \bigg\}.\footnote{If we view $X$ as a matrix, $Z(X)$ is the image of the centered unit cube $C$ under $X$, i.e $X(C)$. }
\end{equation}

When $X$ is full rank, $B_X$ is a compactly supported piece-wise polynomial function. When all the vectors in $X$ lie in the lattice $\Gamma$, the \textit{box spline space associated with $X$} is
\begin{equation}
    \mathcal{S}(X):= \{ B_X \ast' a\, |\, a \in S(\Gamma) \}. 
\end{equation}

Throughout this paper, $X$ will always be full rank and its vectors selected from $\Gamma$. \\

At the foundation of box spline theory is the fact that each $f \in \mathcal{S}(X)$ has two available representations for us to work with. The first, and much more common of the two, is the \textit{$B_X$-representation of $f$}. It is defined as the unique sequence $c_f\in S(\Gamma)$ that satisfies

\begin{equation}
    f = B_X \ast' c_f.
\end{equation}

In order to introduce the second representation of $f\in \cS(X)$, we need to recall that
box spline spaces are intimately connected to specific differential operators and difference operators. Given a non-zero $v \in \R^d$, we define the directional derivative $D_v$ by its action:
\begin{equation}
    D_v := \sum_{i=1}^d v_i \frac{\partial }{\partial x_i}, 
\end{equation}
where $v = (v_1, \dots , v_d)$. More generally, for any multi-set $X$ as above, we define

\begin{equation}
    D_X := \prod_{v \in X} D_v.
\end{equation}

We will consider $D_X$ to act on the space of tempered distributions. We note that applying $D_X$ to any function in $\mathcal{S}(X)$ results in a finite measure (i.e. a sequence) supported on, for our restrictions on $X$ so far, $\Gamma$ or a half shift of $\Gamma$. \\

Next, for any finite collection of vectors $Y \subset \R^d$, define the {\it conical hull of $Y$} as

$${\rm cone}(Y) := \bigg\{\sum_{v \in Y} c_v v\,|\, c_v \geq 0 \ \text{for all} \  v \in Y \bigg\}.$$ 

Suppose that $X \subset \R^d$ is a finite multiset of vectors such that $0$ is the apex of cone$(X)$ \footnote{This can always be achieved given any $X$ by choosing any half-space whose boundary is disjoint of $X$ and flipping the directions of the vectors in $X$ so that all of them lie in the half-space}. The \textit{continuous truncated power associated with $X$} is the distribution $T_X$ such that: 
\begin{enumerate} 
    \item $D_X(T_X) = \delta$. \label{CTP}
    \item The support of $T_X$ lies in cone$(X)$.
\end{enumerate}

For any $f \in \mathcal{S}(X)$, the definition of $T_X$ implies the identity $f=T_X\ast' D_Xf$. This gives us the existence of the \textit{$T_X$-representation of
$f$} defined as the unique sequence $d_f$ such that

\begin{equation}
    f = T_X \ast' d_f
\end{equation}

or equivalently, $d_f := D_Xf$. $T_X$, similar to $B_X$, is a piece-wise polynomial function. However, $T_X \notin \mathcal{S}(X)$ as it's not compactly supported; however, these two representations are connected to each other in a simple way which we describe now. \\

For a non-zero $v \in \R^d$, the \textit{central difference operator} $\Delta_v$ is defined by its action:   
\begin{align}\label{diffopp}
    \Delta_{v}(f) := f(\cdot + \frac{v}{2}) - f(\cdot - \frac{v}{2})
\end{align}
for any $f$ such that \eqref{diffopp} makes sense. More generally, for a multi-set of vectors $X \subset \R^d $, the \textit{central difference operator} $\Delta_X$ is defined as
\begin{align}
    \Delta_X := \prod_{v \in X} \Delta_v .
\end{align}

$\Delta_X$ is a convolution operator and thus we talk about $\Delta_X$ as both a linear map and a sequence in its own right. The following equation holds for any direction set $X$: 

\begin{equation}\label{Derivative_of_Box_Spline}
D_X(B_X) = \Delta_X. 
\end{equation}

Applying $D_X$ to the identity $B_X \ast' c_f = T_X \ast' d_f $ we obtain then that
 
\begin{equation}\label{discreteinnerproduct}
    d_f = \Delta_X c_f.
\end{equation}

We will also assume from here on that each $v \in X$ appears with even multiplicity. We may then consider the following inner product on $\mathcal{S}(X)$ which is at the core of this paper:

\begin{equation}\label{MainInnerProduct}
    \langle f, g \rangle_X :=  (-1)^{|X|/2}\langle c_f , d_g\rangle_{\ell^2(\Gamma)}.
\end{equation}

We study $\cS(X)$ equipped with $\langle \cdot , \cdot \rangle_X $ for particular choices of $X$. In the next section, we discuss how we use $\langle \cdot , \cdot \rangle_X$ to approximate functions from generalized derivative samples.

\subsection{The Approach Taken To Taken In This Paper}
We describe the approach taken in this paper to approximate a smooth and compactly supported $f$ from $\mathcal{S}(X)$ via a discrete collection of generalized derivative samples. We then discuss numerous obstacles with this approach, all of which are overcome later in the paper. \\ 

\textit{Discussion}: In general, let $X_0$ be the multi-set such that $X = [X_0|X_0]$. On $\mathcal{S}(X)$,
\begin{equation}
    \langle D_{X_0}f, D_{X_0}g\rangle_{L^2(\R^d)} 
\end{equation}
defines an inner product. Applying integration by parts gives us that for any $f,g \in \mathcal{S}(X)$,
$$\langle D_{X_0}f, D_{X_0}g\rangle_{L^2(\R^d)} = (-1)^{|X_0|}\langle f|_{\Gamma}, D_Xg \rangle_{\ell^2(\Gamma)} .$$

We consider the multi-set $X_d := [I_d|I_d]$ whose vectors generate the lattice $\Z^d$.\footnote{$B_{X_d}$ is the classic tensor product of linear B-splines i.e $B_{X_d}(x_1, \dots x_d) = \prod_{i=1}^d B_{[1,1]}(x_i)$} Since supp$(B_X) = Z(X)$ for any direction matrix $X$, we have
\begin{equation}
    B_{X_d}|_{\Z^d} = \delta.
\end{equation}

It follows that the restriction $f \mapsto f|_{\Z^d}$ is a bijection from $\mathcal{S}(X_d)$ onto $S(\Z^d)$. Suppose that $\mathcal{B}$ is an (algebraic) orthonormal basis for $\mathcal{S}(X_d)$ under $\langle D_{I_d}\cdot , D_{I_d}\cdot \rangle_{L^2(\R^d)}$. If $f \in \mathcal{S}(X_d)$, then 
\begin{equation}\label{One_Rep}
    f = \sum_{E \in \mathcal{B}} (-1)^d\langle f|_{\Z^d}, D_{X_d}E\rangle_{\ell^2(\Z^d)} E.  
\end{equation}
If $f \in C_c(\R^d)$, the space of continuous and compactly supported functions on $\R^d$, then

\begin{equation}\label{First_Projector}
    P(f) := \sum_{E \in \mathcal{B}} (-1)^d \langle f|_{\Z^d}, D_{X_d}E \rangle_{\ell^2(\Z^d)} E  
\end{equation}

is a projection from $C_c(\R^d)$ onto $\mathcal{S}(X_d)$.  The aforementioned restriction map being bijective and \eqref{One_Rep} imply $\{D_{X_d}E\}_{E \in \mathcal{B}}$ and $\{E|_{\Z^d}\}_{E \in \mathcal{B}}$ form a biorthogonal system, up to a sign, for $S(\Z^d)$. By restricting both sides of \eqref{First_Projector} to $\Z^d$, we obtain that $P$ is the \textit{cardinal interpolation operator} from $C_c(\R^d)$ onto $\mathcal{S}(X_d)$. Equivalently, $P(f)$ is the unique spline in $\mathcal{S}(X_d)$ where
\begin{equation}
    P(f)|_{\Z^d} = f|_{\Z^d}.
\end{equation}

We define the cardinal interpolation operator $Q$ for any $f \in C(\R^d)$ as
\begin{equation}
    Q(f):= B_{X_d} \ast' f|_{\Z^d}.
\end{equation}
$Q(f) = P(f)$ if $f$ is compactly supported, so $Q$ extends the domain of $P$. Approximation properties of $Q$ have been well studied and are typically stated for $f$ lying in

\begin{equation}
    H^m(\R^d) := \{f \in L^2(\R^d) | \ D^{\alpha} f \in L^2(\R^d) \ \text{for all} \  |\alpha| \leq m  \}.
\end{equation}
for some $m > 0$. We endow $H^m(\R^d)$ with the standard norm and semi-norm
\begin{equation}
||f||_{H^m(\R^d)} := \sum_{|\alpha| \leq m} ||D^{\alpha} f||_{L^2(\R^d)}  \ \ \text{and} \  \ |f|_{H^m(\R^d)} := \sum_{|\alpha| = m} ||D^{\alpha} f||_{L^2(\R^d)} .
\end{equation}

For any $h > 0$ and $f \in C(\R^d)$, let $T_h$ be the dilation operator
\begin{equation}
    T_h(f) := f(\cdot/h).
\end{equation}

The following result on cardinal interpolation is classical in one dimension and, for $d > 1$, immediately follows from a result of Carl de Boor (see \cite{BoxSplines}) which we'll state in Section $3$.

\begin{theorem}\label{Schoenberg}
Let $Q_h:= T_hQT_h^{-1}$ and let $m$ be an integer such that $m > d/2$ and $m \geq 2$. For any $f \in H^m(\R^d)$ and $0 < h < 1$, there exists $C$ (depending on $m$ and $d$) such that
\begin{equation}
    ||Q_h(f) - f||_{L^2(\R^d)} \leq Ch^{2}||f||_{H^m(\R^d)}
\end{equation}
\end{theorem}

Now suppose $f$ is $2d$ times continuously differentiable i.e $f \in C_c^{2d}(\R^d)$. For any $A \subset X_d$, integration by parts and Theorem \ref{Schoenberg} imply $P(f)$ approximates $f$ with derivative samples: 
\begin{equation}\label{Derivative_Samples}
    \{\langle D_Af, D_{X_d \backslash A} E\rangle_{L^2(\R^d)}  \}_{E \in \mathcal{B}} \footnote{If $D_{X_d\backslash A} E$ is not a function, we understand $\langle D_Af, D_{X_d \backslash A} E\rangle_{L^2(\R^d)}$ in the distributional sense.}.
    \end{equation}

In this paper, we attempt to use the framework above to approximate smooth compactly supported functions from generalized derivative samples as in \ref{Derivative_Samples} for much more general $X$. In particular, we attempt to use two algebraic bases $\mathcal{B} = \{E\}_{E \in \mathcal{B}}$ and $ \mathcal{B}'= \{\tilde{E}\}_{E \in \mathcal{B}}$ for $\mathcal{S}(X)$, construct the analysis-synthesis operator

\begin{equation}\label{Straight_Forward_Map}
    f \mapsto \sum_{E \in \mathcal{B}} (-1)^{|X_0|} \langle f|_{\Gamma}, D_XE \rangle_{\ell^2(\Gamma)} \tilde{E},  
\end{equation}

and then convert the coefficients above into generalized derivative samples. There are several major obstacles with this approach which we now describe. \\

First and foremost, given an arbitrary multi-set $X$ and bases $\mathcal{B}$ and $\mathcal{B}'$ for $\mathcal{S}(X)$, the operator \eqref{Straight_Forward_Map} may not be well-defined for all $f \in C_c(\R^d)$ (even if a basis orthonormal under $\langle D_{X_0}\cdot, D_{X_0}\cdot \rangle_{L^2(\R^d)}$ is used). Secondly, if \ref{Straight_Forward_Map} is well-defined, it is no longer immediately clear if it approximates $f$. Third, there is one more important issue here which underlies any attempt at using generalized derivative samples of a function $f$ to approximate it in $L^2$ norm. \\

 If $\mathcal{B}$ is a collection of functions for which $\eqref{Straight_Forward_Map}$ approximates $f$ via the approximation scheme in Theorem \ref{Schoenberg}, and the coefficients can be changed to derivative samples as in \eqref{Derivative_Samples}, then the set of diameters of $\{\text{supp}(D_AE)\}_{E \in \mathcal{B} }$ cannot be uniformly bounded for any $A \subsetneq X$. This is one of the main reasons that most work on approximations with derivative samples use function samples in combination (e.g \cite{GarciaPerezVillalon2008}, \cite{UnserZerubia1997}, \cite{UnserZerubia1998}, \cite{PriyankaSelvan2024}, \cite{Corso2025}). Thus, if the first two issues above are overcome, the inherent non-locality of $\mathcal{B}$ means its functions or the coefficients of the approximation map may be computationally undesirable.  \\

Our resolution to these obstacles is as follows. First, we seek to construct an orthonormal basis $\mathcal{B}$ for $\mathcal{S}(X)$ under $\langle \cdot , \cdot \rangle_X$ (more generally we use a biorthogonal system). $\langle \cdot, \cdot \rangle_X$ is similar to $\langle D_{X_0}\cdot, D_{X_0}\cdot \rangle_{L^2(\R^d)}$\footnote{The two inner products are equivalent whenever $B_X|_{\Gamma}$ is invertible as a convolution operator.} and results in \eqref{Straight_Forward_Map} being a well-defined finite sum for all $f \in C_c(\R^d)$. Secondly, we look for all $E \in \mathcal{B}$ to satisfy $|$supp$(D_XE)| < N$ for some positive integer $N$. This results in each $E \in \mathcal{B}$ being, at most, a linear combination of $N$ shifts of $T_X$. Thus, $\mathcal{B}$ would be local in its $T_X$-representation and the outlined approach may be worthy of pursuit.\\

However, finding an orthogonal basis under $\langle \cdot , \cdot \rangle_X$ with uniformly finite $T_X$-representations for a large collection of spline spaces $\mathcal{S}(X)$ seems, if even possible, highly non-trivial. This is what makes the next section remarkable and the foundational motivation of this paper.

\subsection{\hspace{-8mm} . Weyl's Character Formula - Box Spline Version}

Suppose that the support of $f \in \mathcal{S}(X)$ is a compact convex polyhedron $K$. Trivially, supp$(d_f)$ includes the extreme points of $K$. Moreover, it seems exceptionally difficult to find examples of $f \in \mathcal{S}(X)$ where supp$(d_f)$ consists {\it only} of extreme points of $K$. However, in the realm of non-commutative algebra lies a theory which provides us with a wealth of such examples for particular choices of $X$. This theory is rooted in the celebrated \textit{Weyl's Character Formula}, which we re-formulate in the language of box splines below.  \\

Throughout this paper $R \subset \R^d$ is a \textit{root system} (\cite{Hall}) and we assume:

\begin{itemize}
    \item $W$ is the Weyl group of $R$ (i.e., the group of reflections generated by $R$) 
    \item $R^+$ is a set of positive roots in $R$
    \item $\Pi$ is the set of simple roots for $R^+$.  
    \item $\Lambda$ is the root lattice.
    \item $C$ is the (open) \textit{fundamental Weyl chamber} associated to $R$ i.e, 
    $$C :=\{v \in \R^d |   v \cdot r > 0 \ \text{for all} \  r \in R^+ \}. $$
    
    \item  For any $v \in \R^d$, $v^*$ is the unique element in $W(v)$ (the orbit of $v$ under $W$) that lies in the closure of $\overline{C}$ of $C$. More generally $P^*:=\{v^*| v \in P \}$ for any $P \subset \Lambda$. 
    \item  $\geq$ is the following ordering on $\Lambda^*$:  $\lambda \geq \mu$ if and only if conv$(W(\mu)) \subseteq $ conv$(W(\lambda))$. This is equivalently to the statement that $\lambda^* - \mu^* \in \text{cone}(R^+)$.
\end{itemize}

Throughout this paper, we often work with functions which are supported on the convex hull of $W(\lambda)$ where $\lambda \in \Lambda$. Such a polytope is often referred to as a \textit{Weyl polytope} and $\Lambda^*$ serves as a choice of representatives for the set of all possible Weyl polytopes.  \\

We recall that each $w\in W$ has an associated {\it length} $l(w) \in \Z_{+}$ and we define the \textit{Weyl vector}

\begin{equation}
    \overline{\rho} := \frac{1}{2} \sum_{r \in R^+} r
\end{equation}

\begin{theorem}\label{splineversionofweyl}
Let $\lambda \in \Lambda^*$. The function

\begin{equation}\label{SymmetricWeylSpline}
S_\lambda:=\sum_{w\in W} (-1)^{l(w)} T_{R^+}(\cdot-w(\lambda + \overline{\rho}))
\end{equation}

is compactly supported. Thus, for each such $\lambda$, the spline $S_\lambda$ has the following two properties:
\begin{itemize}
    \item It is supported in the convex hull of $W(\lambda + \overline{\rho})$.
    \item It satisfies: supp$(D_{R^+}S_\lambda) = W(\lambda + \overline{\rho})$.
\end{itemize} 
\end{theorem}

Weyl's character formula, in its original form, has a significantly different context than the above spline version. Relevant to the original form of the statement is that if
\begin{equation}
S_\lambda=B_{R^+}\ast' s_\lambda:=\sum_{\mu\in\Lambda}s_\lambda(\mu)B_{R^+}(\cdot-\mu),
\end{equation}
the value $s_\lambda(\mu)$ is the multiplicity of the eigenvector $\mu$ in a suitable irreducible representation of the semi-simple Lie algebra associated with $R$. Similarly, $S_{\lambda}$, under an appropriate isomorphism, is essentially the \textit{Schur polynomial} corresponding to $\lambda$. This group  representation connection, while relevant, is not truly pertinent to the content of this paper. \\

\smallskip{\it Examples:} Up to a possible linear transformation, the standard 2D 3-directional and 4-directional box splines are associated with a root system, hence are covered by Theorem \ref{splineversionofweyl} (and by the actual results of this paper, as well!).
\begin{itemize}
    \item For the (piecewise linear) 3-directional box spline space in 2D, the Weyl group is of order 6, and so there exists compactly supported splines whose supports have arbitrarily large diameter and are a linear combination of shifts of 6 truncated powers.
    \item For the (piecewise quadratic) 4-directional box spline space in 2D, the Weyl group is of order 8, and so there exists compactly supported splines whose supports have arbitrarily large diameter and are a linear combination of shifts of 8 truncated powers.
\end{itemize}

\smallskip{ \it Orthogonality of the above Weyl's character splines.}  \\ 

Let us equip the space span$\{S_{\lambda}\}_{\lambda \in \Lambda^*}$ with the following scalar product:  

\begin{equation}\label{Symmetric IP}
(S_\lambda,S_\mu)_{R^+} := \langle D_{R^+}S_{\lambda}, D_{R^+} S_{\mu}  \rangle.
\end{equation}

If $x,y \in C$ are distinct, the orbits $ W(x)$ and $W(y)$ are disjoint. It follows that

\begin{equation}\label{SymmetricSplineOrthogonality}
( S_\lambda,S_\mu)_{R^+}= 0.
\end{equation}

The orthogonal collection $\{S_{\lambda}\}_{\lambda \in \Lambda^\ast}$ falls short of being an orthogonal basis for the entire spline space $\mathcal{S}(R^+)$ as each $S_{\lambda}$ is invariant under the action of $W$. However, it is an orthogonal basis for the $W$-invariant subspace $\mathcal{S}(R^+)_W$ of $\mathcal{S}(R^+)$. Moreover, it is easy to see that 

\begin{equation}
    (S_{\lambda}, S_{\mu})_{R^+} = (-1)^{|R^+|} \langle \Delta_{R^+}^2 s_{\lambda}, s_{\mu} \rangle  
\end{equation}
and the right hand side bears similarity to $\langle \cdot, \cdot \rangle_X $ in the case where $X = [R^+| R^+]$. \\

It is these observations mentioned in this section that motivated the pursuit of one of the achievements of this paper: An orthogonal basis for $\mathcal{S}(X)$ under $\langle \cdot, \cdot \rangle_X$, with each element having sparse $T_X$-representation, for a variety of choices of $X$, all related to root systems.

\subsection{\hspace{-8mm} . The Orthogonal Splines In This Paper}

Since the proof of the Weyl Character Formula, significant progress has been made in the study of characters of irreducible representations. As a consequence, a plethora of generalizations of the Schur polynomials have been established and have found uses in a variety of fields. These polynomials, when equipped with a suitable scalar product, bear orthogonality relations which form the foundation for a general theory of orthogonal polynomials associated to root systems which is described in \cite{Mac00}.  \\

We draw on this theory, particularly the Jack polynomials \cite{Jack69} and Heckman-Opdam polynomials (\cite{HeckmanI} and \cite{HeckmanII}) in a similar way that we drew upon the Weyl character formula. More specifically, we use this theory to introduce a collection of inner products on particular box spline spaces which bear far more powerful and general orthogonality properties than what we witnessed with $( \cdot, \cdot )_{R^+}$ on $\mathcal{S}(R^+)_W$. To describe these box spline spaces, we introduce the following notation which will be consistent throughout the rest of this paper. \\

We denote by $X_0\subset \R^d$ a fixed multi-set that satisfies the following two conditions: 
\begin{enumerate}
    \item Every positive root $r \in R^+$ lies in $X_0$ and appears with multiplicity $k_{r} > 0$. 
    \item Positive roots in the same orbit under $W$ have the same multiplicity in $X_0$. 
\end{enumerate}

We then define another multi-set $X$ to be composed of two copies of $X_0$ i.e 
\begin{equation}
X:= [X_0|X_0].
\end{equation}
The splines we study in this paper are those lying in $\mathcal{S}(X)$. Before describing the aforementioned orthogonality property held by $\mathcal{S}(X)$ under $\langle \cdot, \cdot \rangle_X$, we give a few remarks. For these choices of $X$, the inner products $\langle \cdot , \cdot \rangle_X$ are a sub-family of a much larger family of scalar products on certain \textit{exponential box spline spaces} rooted in the work of Macdonald (see \cite{Mac96} and \cite{Mac00}), but we do not discuss this in the present paper. Others have considered box splines associated with root systems previously (see \cite{Hyperplane}, \cite{mcSwiggen2021}, and \cite{Minho11}), but to the knowledge of the author, this particular connection with orthogonal polynomials has not been explored before. \\

In \cite{O}, a remarkable orthogonal basis of polynomials is constructed with properties which imply, via isomorphism, an orthogonal decomposition of $\mathcal{S}(X)$ into spaces determined by the orbits of $\Lambda$ under $W$. More specifically, we have an orthogonal decomposition under $\langle \cdot , \cdot \rangle_X$

\begin{equation}\label{ODecomposition}
    \mathcal{S}(X) = \bigoplus_{\lambda \in \Lambda^*} \mathcal{S}_X(\lambda),  
\end{equation}

where each $\mathcal{S}_X(\lambda)$ is $W$-invariant, has dimension $|W(\lambda)|$, and only contains splines supported on a polytope slightly larger than $\text{conv}(W(\lambda)) $. We refer to the splines in  $\mathcal{S}_X(\lambda)$ as \textit{Weyl splines} and define them precisely later on. In the appendix, we discuss the theory in \cite{O} and how it implies the orthogonal decomposition above. \\

Earlier, we saw that the sparsity and shape of the $T_{R^+}$-representation of $\{S_{\lambda} \}_{\lambda \in \Lambda^*}$ implied an orthogonality relation. In a reversal to this, we'll find that the orthogonality of \eqref{ODecomposition} forces any basis respecting the decomposition to have a sparse $T_X$-representation.

\subsection{Goals and Outline of this paper}
Throughout the rest of this paper, we assume the notation related to root systems of Section 1.3. The goal of this paper is to investigate bases of $\mathcal{S}(X)$ which respect the orthogonal decomposition in \eqref{ODecomposition} and use these bases to construct maps similar to \eqref{Straight_Forward_Map} from which we can approximate smooth $f \in C_c(\R^d)$ from derivative samples induced by these bases.\\

In Section 2, we introduce the Weyl spline spaces $\{\mathcal{S}_X(\lambda)\}_{\lambda \in \Lambda^*}$ which comprise the decomposition \eqref{ODecomposition}. It'll be clear that the splines in $\mathcal{S}_{X}(\lambda)$ are not local in the usual sense as the diameter of supp$(E)$ for any $E \in \mathcal{S}_X(\lambda)$ is proportional to the distance of $\lambda$ from the origin. The main result of Section $2$ is that these splines \textit{are} local in a sense which we describe now. For any finite multi-set $S \subset \Lambda$, define the vector 

\begin{equation}
    \rho_S := \sum_{s \in S} s. 
\end{equation}

For any generic choice of $\lambda \in \Lambda^*$ and $E \in \cS_X(\lambda)$, 

\begin{equation}\label{LocalT_XRep}
{\rm supp}(D_XE) \subset \bigcup_{w \in W} w\bigg (\lambda + \frac{\rho_{X_0}}{2} + Z(X_0) \bigg).
\end{equation}

So, the size of the support of $D_XE$ is determined by $X_0$ and $W$ with $\lambda$ only influencing the positioning of supp$(D_XE)$. As $E = T_X \ast' D_XE$, this implies that $E$ is a linear combination of a relatively small number of shifts of $T_X$. \\

In Section 3, we take a basis of Weyl splines $\mathcal{B}$ for $\mathcal{S}(X)$ and its dual basis $\mathcal{B}^*$ under $\langle \cdot, \cdot \rangle_X$ and construct approximation maps into $\mathcal{S}(X)$. For $a \in S(\Lambda)$, we consider $P_a:C_c(\R^d) \to \mathcal{S}(X)$,

\begin{equation}\label{Quasi_Interpolant}
    P_a(f) := \sum_{E \in \mathcal{B}} (-1)^{|X_0|}\langle a \ast f|_{\Lambda}, D_XE \rangle_{\ell^2(\Lambda)} E^*. 
\end{equation}
The operator is well-defined and, we show, is a local and discrete \textit{quasi-interpolation operator} on $\mathcal{S}(X)$. Moreover, if the root system underling $X$ is irreducible, $a \in S(\Lambda)$ can be chosen so that $P_a$ has optimal approximation order on sufficiently smooth functions. Secondly, we show that if $X$ comes from the $A_n$ root system, we can extend the domain of $\langle \cdot, \cdot \rangle_X$ to a larger space where it becomes a bilinear form $(\cdot, \cdot)_X$. If $f$ is compactly supported, the map
\begin{equation}\label{First_projector}
    \mathcal{T}: f \mapsto \sum_{E \in \mathcal{B}} (f, E)_XE^*
\end{equation}
is a projector and also a quasi-interpolation operator with optimal approximation order. Arguably most important is that both $P_a$ and $\mathcal{T}$, after converting the coefficients of \eqref{Quasi_Interpolant} and \eqref{First_projector} with integration by parts, can be expressed in terms of generalized derivative samples:  
\begin{equation}
    \{\langle D_Af, D_{X \backslash A} E \rangle_{L^2(\R^d)}\}_{E \in \mathcal{B}}
\end{equation} 

where $A$ varies over a finite collection of multi-sets. \\

In Section $4$, we investigate these bases more computationally. In particular, we detail how to explicitly construct biorthogonal systems for $\mathcal{S}(X)$ composed of Weyl splines. This first requires us to compute $\langle \cdot, \cdot \rangle_X$ on $\cS_X(\lambda)$. The work in Section $2$ implies that for any $f,g \in \mathcal{S}_X(\lambda)$, supp$(c_f)$ and supp$(d_g)$ intersect only on $W(\lambda)$ and so computing $\langle f , g\rangle_X$ can be done by computing $D_Xf$ at $W(\lambda)$. We give a formula for these values, for a generic choice of $\lambda$, and use this formula to compute a biorthogonal systems of Weyl splines. We'll see that some systems are more desirable than others. Examples are given in one dimension.

\section{$\cS_X(\lambda)$ and The Weyl Spectrum: The Support of $D_X\cS_X(\lambda)$}
For the entirety of this section, we fix $\lambda \in \Lambda^*$, $X_0$ as in Section 1.4, $X = [X_0|X_0]$, and
\begin{equation}
    \rho := \rho_{X_0}. 
\end{equation}
\def\lam{\lambda}

We define $\mathcal{S}_X(\lambda)$ as the set of all $E\in \mathcal{S}(X)$ satisfying the following:

\begin{enumerate}
    \item  $\text{supp}(E) \subseteq \text{conv}(W(\lambda + \rho)).$ \label{Property 1}
    \item  $(D_XE) = 0$ on $\text{conv}(W(\lambda)) \backslash W(\lambda)$. \label{property 2} 
\end{enumerate}

We refer to any element of $\mathcal{S}_X(\lambda)$ as a \textit{Weyl spline} associated to $X$ and $\lambda$. It will be very useful to understand the Weyl splines in terms of their $B_X$ representation. Recall that for two sets $A,B \subset \R^d$, their Minkowski sum is defined as: 
\begin{equation}
    A+B := \{a+b\ | \ a \in A \ \text{and} \ b \in B  \}.
\end{equation}

It is well known that the Minkowski sum of two convex sets is also convex. Moreover, if $A$ and $B$ are invariant under $W$, so is $A+B$.  It follows that

\begin{equation}
   \text{conv} ( W(\lambda + \rho) ) \subseteq \text{conv} ( W(\lambda)) + \text{conv}(W(\rho)).
\end{equation}


Recall that $\text{conv}(A) + \text{conv}(B) = \text{conv}(A+B)$ for any $A,B \subset \R^d$. It is known that $\rho \in \Lambda^*$ and so $\lambda + \rho \in \Lambda^*$. It follows that $\mu \in W(\lambda) + W(\rho)$ satisfies $\mu^* \leq \lambda + \rho$ and thus lies in conv$(W(\lambda + \rho))$ by Lemma \ref{lem:convdom}. This implies  
\begin{equation}
   \text{conv} ( W(\lambda)) + \text{conv}(W(\rho) ) \subseteq \text{conv} ( W(\lambda + \rho) ).
\end{equation}

Since supp$(B_X)$ = $Z(X)$ and $Z(X) = \text{conv}(W(\rho))$, we may again use Lemma \ref{lem:convdom} to conclude that Property 1 in the definition of $\mathcal{S}_X(\lambda)$ is equivalent to the statement:

\begin{equation}\label{SplineSupport}
    E \in \text{span}\{B_X(\cdot-\mu) \}_{\mu^* \leq \lambda}. 
\end{equation}

We may also equivalently restate property \ref{property 2} in the definition of $\mathcal{S}_X(\lambda)$ with the condition: 
\begin{equation}\label{Property 2 Restatement}
D_X(E) = 0 \quad {\rm on}\quad \{\mu \in \Lambda| \ \mu^* < \lambda \}.
\end{equation}

\begin{figure}[h!]
\centering
\includegraphics[scale =.53]{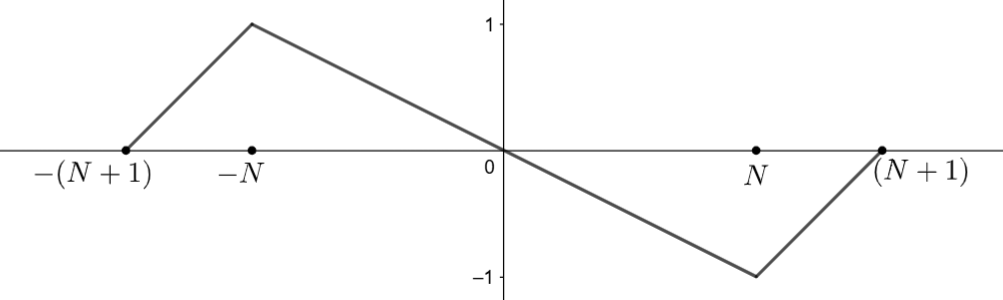}
\caption{The figure above is an example in 1D of a spline $S \in \mathcal{S}_X(N)$ where $X = [1,1]$ and $N$ is any (fixed) positive integer. In this case, the Weyl group $W$ consists of the identity and the map $x \mapsto -x$. The sparsity of $d_S$ is clear here as supp$(D^2S) = \{\pm N,\pm (N+1) \}$. }
\end{figure}

It is not hard to show that $\mathcal{S}_X(\lambda)$ is a linear subspace of $\mathcal{S}(X)$ of dimension $|W(\lambda)|$ and we will prove a slightly stronger statement later. It is also not hard to show that
                    
\begin{equation}\label{eq: Decomposition}
    \mathcal{S}(X) = \bigoplus_{\lambda \in \Lambda^*} \mathcal{S}_X(\lambda).
\end{equation}   

On the other hand, the following result, implied by the work done in \cite{O}, is highly non-trivial, paramount throughout this paper, and we give its details in the appendix. 
\begin{theorem}
    The following direct sum decomposition is orthogonal under $\langle \cdot, \cdot \rangle_X$:
    \begin{equation}
    \mathcal{S}(X) = \bigoplus_{\lambda \in \Lambda^*} \mathcal{S}_X(\lambda).
    \end{equation}
\end{theorem}

For any $\lambda \in \Lambda^*$, we define the following orbit of zonotopes
\begin{equation}
    Z_{X_0}(\lambda) := \bigcup_{w \in W} w\bigg(\lambda + \frac{\rho}{2} + Z(X_0) \bigg).
\end{equation}

We define the \textit{Weyl Spectrum} $W_{X_0}(\lambda)$ of $\mathcal{S}_X(\lambda)$ as the discrete points of these zonotopes: 

\begin{equation}
W_{X_0}(\lambda) := \Lambda \cap Z_{X_0}(\lambda).
\end{equation}

It is easy to show and helpful to know that 
\begin{equation}
    W_{X_0}(\lambda) = \{w(\lambda + \rho_A) | A \subset X_0 \  \text{and} \ w \in W\}.   
\end{equation}

With this, we now state the main result of this section.

\begin{theorem}\label{GeneralizedWeylCharacter}
Let $\lambda \in \Lambda^*$ satisfy that if $\lambda \leq \mu \leq \lambda + \rho$, then $\mu^*$ lies within the open fundamental Weyl chamber $C$. Then for any $E \in \mathcal{S}_X(\lambda)$, ${\rm supp}(D_X E )\subseteq W_{X_0}(\lambda)$.\footnote{Note that, unlike conv$(W(\lambda + \rho))$, $W_{X_0}(\lambda)$ is sparse:  The number of lattice points in the latter depends only on $X_0$, hence stays fixed as $\lambda$ increases.} 
\end{theorem}

The condition on $\lambda$ that $\lambda \leq \mu  \leq \lambda + \rho $, $\mu^* \in C$ is equivalent to $\lambda$ being a fixed distance away from the boundary of $C$. More precisely, given a ball of radius $r$ centered at the origin, the number of $\lambda \in \Lambda^*$ contained in the ball is $O(r^d)$ while the number of points excluded by the condition is $O(r^{d-1})$. So, the theorem holds for a generic choice of $\lambda \in \Lambda^*$. In proving this theorem, we will still see a weaker uniform bound on $|$supp$(D_XE)|$ for all $\lambda \in \Lambda^*$. 

\begin{figure}[h!] \label{Fig2}
\begin{subfigure}[h!]{0.4\linewidth}
\includegraphics[width = 5.5cm ,height = 5.4cm]{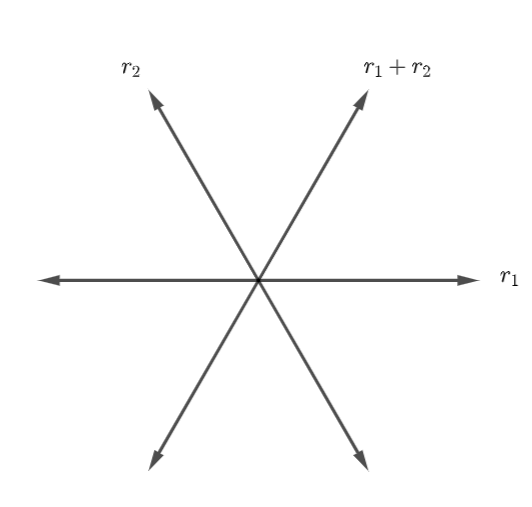}
\includegraphics[scale = .110]{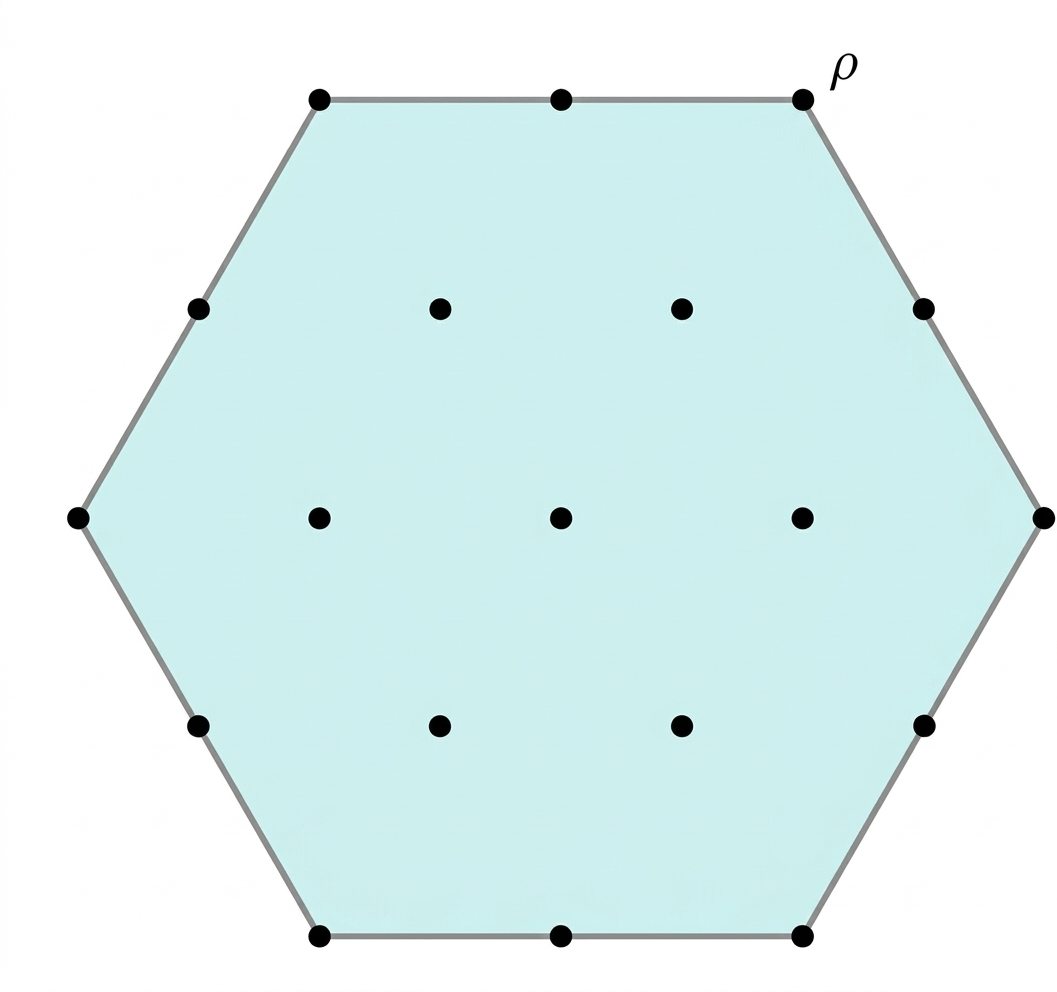}
\end{subfigure}
\hspace{-4mm}
\begin{subfigure}[h!]{0.4\linewidth}
\includegraphics[scale =.80]{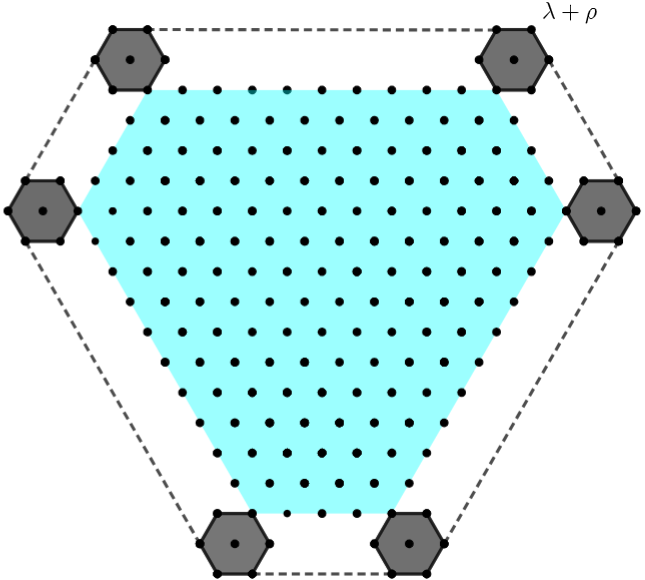}
\end{subfigure}

\caption{The top left part shows the $A_2$ root system with the three positive roots labeled. The bottom left shows $Z(X)$, the support of the box spline $B_X$, where each positive root in $X$ appears with multiplicity $2$ (so $|X| = 6$). On the right, the region within the dotted lines represents the maximal support of a generic $E \in \mathcal{S}_X(\lambda)$ and the lightly blue shaded region is the corresponding maximal support of $c_E$. The darker shaded hexagons represent the Weyl spectrum of $E$, i.e., the maximal support of $d_E$.  } 
\label{Fig2}
\end{figure}

\newpage
Since $D_X$ maps $\mathcal{S}(X)$ to $S(\Lambda)$, we know from \eqref{SplineSupport} that supp$(D_XE)$ is contained in the set
\begin{align}\label{originalsupp}
  \{\mu \in \Lambda| \  \mu^* \leq \lambda + \rho \} 
\end{align}

To prove Theorem \ref{GeneralizedWeylCharacter}, we start with the set in \eqref{originalsupp} and provide a series of arguments, each of which subsequently excludes points from it, in a way that the remaining set still contains supp$(D_XE)$. The final set we obtain in this sequence of exclusions is $W_{X_0}(\lambda)$. \\
 
By Property $\ref{property 2}$ in the definition of $\cS_X(\lambda)$, $D_XE = 0$ on conv$(W(\lambda))\backslash W(\lambda)$. Using \eqref{Property 2 Restatement} and removing the appropriate points from \eqref{originalsupp}, supp$(D_XE)$ lies in the following discrete annulus:  
\begin{align}\label{annlus}
A_{X_0}(\lambda) := \{ \mu \in \Lambda \ |\  \mu^*  \not< \lambda \ \text{and} \ \mu^* \leq \lambda + \rho \}.
\end{align}
The following lemma narrows down the support of $D_XE$ from $A_{X_0}(\lambda)$ to a significantly smaller set and, in particular, implies that $|$supp$(D_XE)|$ is independent of the length of $\lambda$. 
\begin{lemma}\label{SignificantReduction} 
   If $E \in \mathcal{S}_X(\lambda)$, then supp$(D_XE)$ is contained in the set 
    \begin{align}\label{Parallelogram}
    P_{X_0}(\lambda) :=
    \{\mu \in \Lambda \  |\ \lambda  \leq \mu^* \leq \lambda + \rho  \}.
\end{align}
\end{lemma}

\begin{proof}

Let $\mu \in A_{X_0}(\lambda)\backslash P_{X_0}(\lambda)$. We will show that $D_XE(\mu) = 0$. Define    
\begin{align}
    V_{\mu^*} &:= \text{span}\{B_X(\cdot - \gamma) \}_{\gamma^* \leq \mu^*}.
\end{align}

Property \ref{Property 1} in the definition of the Weyl spline space implies that if $v \in \Lambda^*$ and $v \leq \mu^*$, then $\mathcal{S}_X(v) \subset V_{\mu^*}$. Since the spaces in $\{\mathcal{S}_X(v)\}_{v \in \Lambda^*}$ are orthogonal, they are certainly linearly independent as subspaces from one another. Of course, the box splines $\{B_X(\cdot - \gamma) \}_{\gamma^* \leq \mu^*}$ are also linearly independent, and so dimension considerations imply that

\begin{equation}
    V_{\mu^*} = \bigoplus_{\{v \leq \mu^*| \ v \in \Lambda^* \}} \mathcal{S}_X(v) .  
\end{equation}

Since $\lambda \not\leq \mu^*$, $\mathcal{S}_X(\lambda)$ is orthogonal to $V_{\mu^*}$. In particular, $B_X(\cdot - \mu) \in V_{\mu^*}$, so we have  
\begin{align}
    0 & = \langle E, B_X(\cdot - \mu) \rangle_X \\
      & = (-1)^{|X_0|}D_XE(\mu)    
\end{align}
which completes the proof.
\end{proof}

\begin{figure}[h!]
\caption{Below are the three sets of points: \eqref{originalsupp}, $A_{X_0}(\lambda)$, and $P_{X_0}(\lambda)$ for the specific case where the root system is $B_2 \subset \R^2$, the positive roots are $R^+ =\{(1,0),(1,1), (0,1), (-1,1) \}$, and $X_0 = R^+$. $\lambda \in \Lambda^*$ satisfies the conditions of Theorem \ref{GeneralizedWeylCharacter}, but is arbitrary otherwise.}

\begin{subfigure}[h!]{0.4\linewidth}
\hspace{.5cm}
\includegraphics[width = 6.1cm ,height = 5cm]{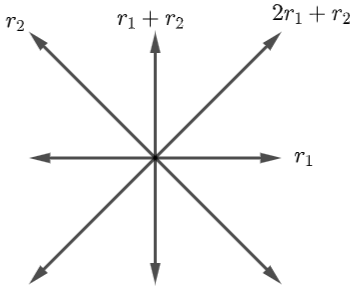}
\vspace{3mm}
\caption*{(a) The $B_2$ root system with our choice of positive roots labeled.}
\end{subfigure}
\hspace{2cm}
\begin{subfigure}[h!]{0.4\linewidth}
\vspace{6mm}
\includegraphics[width = 8cm ,height = 8cm, trim= 0cm 0cm 0cm 5mm]{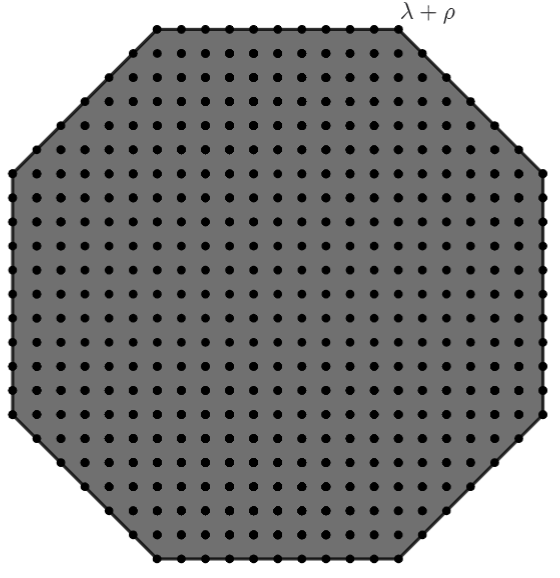}
\caption*{(b) The first set containing supp$(D_XE)$ as stated in \eqref{originalsupp}. }
\end{subfigure} \\
\end{figure}

\begin{figure}[h!]
\begin{subfigure}[h!]{0.4\linewidth}
\hspace{-1cm}
\includegraphics[width = 8cm ,height = 8.2cm]{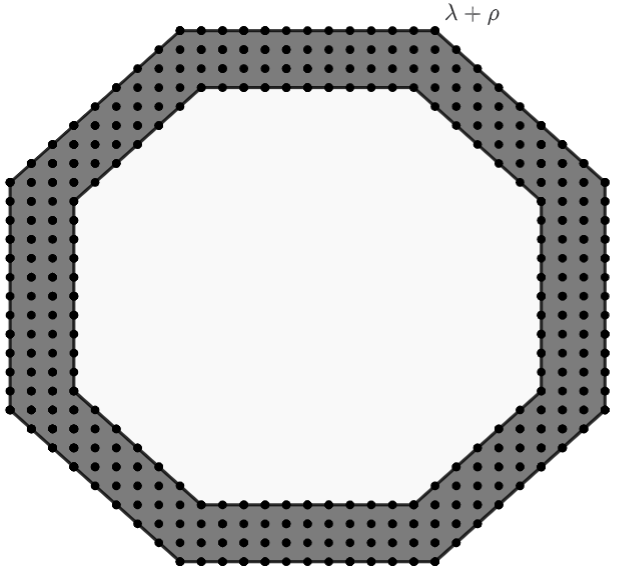} 
\caption*{(c) The annulus $A_{X_0}(\lambda)$ obtained from excluding points from \eqref{originalsupp}.}
\end{subfigure}
\hspace{2cm}
\begin{subfigure}[h!]{0.4\linewidth}

\includegraphics[width = 8cm ,height = 7.7cm]{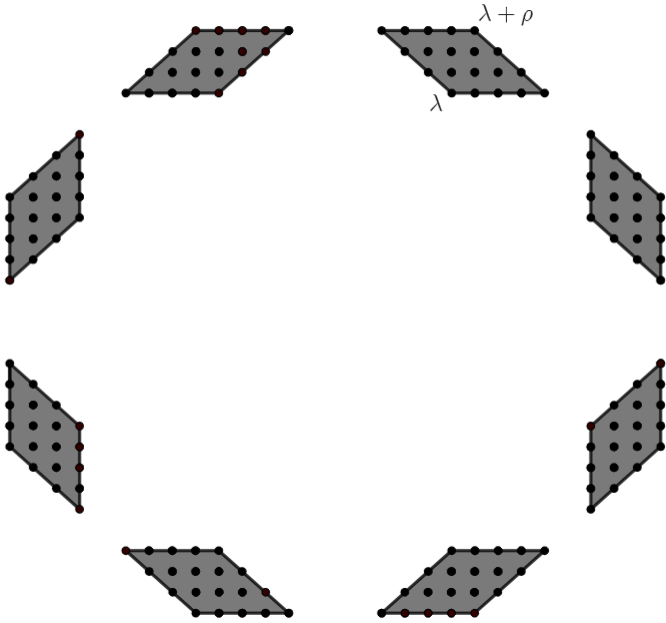}
\caption*{(d) \ $P_{X_0}(\lambda)$ which is obtained by removing points from $A_{X_0}(\lambda)$.}
\end{subfigure} 
\end{figure}

\newpage

Now we work towards using that supp$(D_{X}E) \subset P_{X_0}(\lambda)$ to show supp$(D_XE)$ is contained in the smaller set $W_{X_0}(\lambda)$. To this end, we make use of properties of $T_{X}$ as follows. For each $r \in R^+$, let $X_{r}$ be the multi-set obtained by removing all vectors from $X$ except $r$ and let 
\begin{equation}
E_{r} := D_{X \backslash X_{r}}E.
\end{equation}


By expressing $E$ in terms of its $T_X$-representation and using the fact that $D_YT_X = T_{X \backslash Y}$ for any $Y \subset X$, we can use that supp$(D_XE) \subset P_{X_0}(\lambda)$ to conclude 

\begin{equation}
E_{r} = \sum_{\mu \in P_{X_0}(\lambda)} d_E(\mu) T_{X_r} (\cdot - \mu). 
\end{equation}

We refer to a line parallel to some $r \in R$ as a \textit{root line}. Suppose $L_{r}$ is a root line parallel to $r \in R$ and having non-empty intersection with $P_{X_0}(\lambda)$. Moreover, recall that the multiplicity of $r$ in $X$ is $2k_r$. Since $E$ is compactly supported, $E_r$ is compactly supported. It follows that 
\begin{equation}
    \sum_{\mu \in P_{X_0}(\lambda) \cap L_r} d_E(\mu) T_{X_r} (\cdot - \mu)
\end{equation}

is also compactly supported. Lemma \ref{Truncated_Power_Lemma} implies that if $|L_{r}\cap P_{X_0}(\lambda)| < 2k_{r}+1$, then $d_E = 0$ on $L_r \cap P_{X_0}(\lambda)$. This observation will be critical for narrowing down the containment of supp$(D_XE)$ from $P_{X_0}(\lambda)$ to $W_{X_0}(\lambda)$. Now, we fix 

\begin{equation}
    P_{\lambda}:=  \{\mu \in \Lambda^*| \ \lambda  \leq \mu \leq \lambda + \rho  \}.
\end{equation}

The assumption made in Theorem \ref{GeneralizedWeylCharacter} is that the above $P_{\lambda}$ lies in the open fundamental Weyl chamber $C$. We now adopt this assumption from here on. From this assumption it follows that the sets in $W(P_{\lambda})$, the orbit of $P_{\lambda}$ under $W$, are pairwise disjoint. \\ 

In this case, any root line intersecting $P_{X_0}(\lambda)$ does so at an even number of points. Thus, the existence of a root line $L_{r}$ such that $|P_{X_0}(\lambda) \cap L_{r}| < 2k_{r}+2$ is equivalent to the existence of $w \in W$ and a root line $L_{w(r)}$ intersecting $P_{\lambda}$ such that $|L_{w(r)} \cap P_{\lambda}| <k_{w(r)}+1$.  In lieu of this, we let $L_{r}(\mu)$ be the line parallel to $r$ and containing $\mu$ and define for any $S \subset \Lambda$
\begin{equation}
    L(S):= \{\mu \in S | \  |L_{r}(\mu) \cap S| \geq k_{r} +1\ \text{for all} \ r \in R \}.
\end{equation}

Clearly $L(S) \subset S$ for any $S \subset \Lambda$. We may recursively apply $L$ to $P_{\lambda}$ and conclude from our observations above that for any positive integer $i$, 
\begin{equation}
    \text{supp}(D_XE) \subset \bigcup_{w \in W} w (L^i(P_{\lambda})).
\end{equation}
Since $P_{\lambda}$ is a finite set, there exists a smallest integer $m$ such that
\begin{equation}
L^m(P_{\lambda}) = L^{m+1}(P_{\lambda}). 
\end{equation}
We prove Theorem \ref{GeneralizedWeylCharacter} by showing $L^m(P_{\lambda}) = W_{X_0}(\lambda)^*$. We first require two lemmas.

\begin{figure}[h!]
\caption{Below, we consider the same example that was done in Figure $3.3$ and show, for this particular case, a visualization for how we narrow down the support of $D_XE$ from $P_{X_0}(\lambda)$ to $W_{X_0}(\lambda)$. This ultimately proves Theorem \ref{GeneralizedWeylCharacter} in a specific case. }

\begin{subfigure}[h!]{0.4\linewidth}
\hspace{-1cm}
\vspace{.6cm}
\includegraphics[width = 8.1cm ,height = 8.3cm]{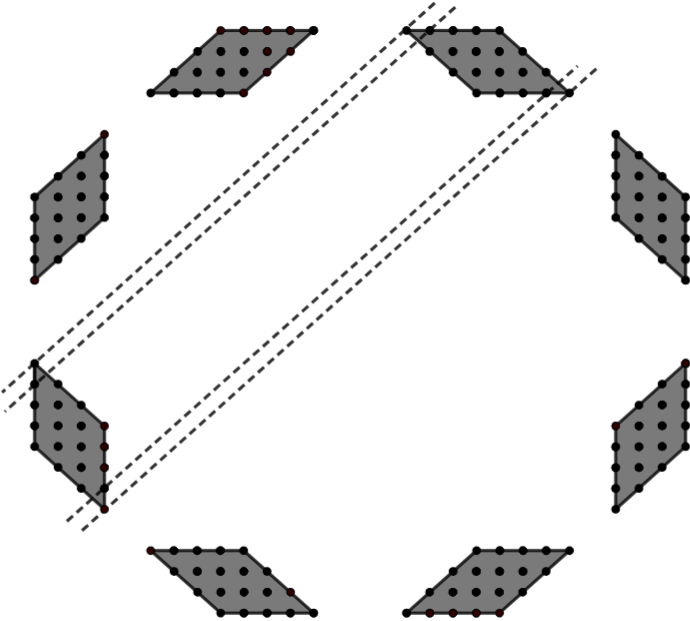}
\vspace{-1cm}
\caption*{(a) $P_{X_0}(\lambda) = W(P_{\lambda})$ with root lines containing points not in $L(P_{\lambda})$.   }
\end{subfigure}
\hspace{3cm}
\begin{subfigure}[h!]{0.4\linewidth}
\hspace{-1cm}
\includegraphics[width = 8cm ,height = 8cm, trim= 0cm 0cm 0cm 5mm]{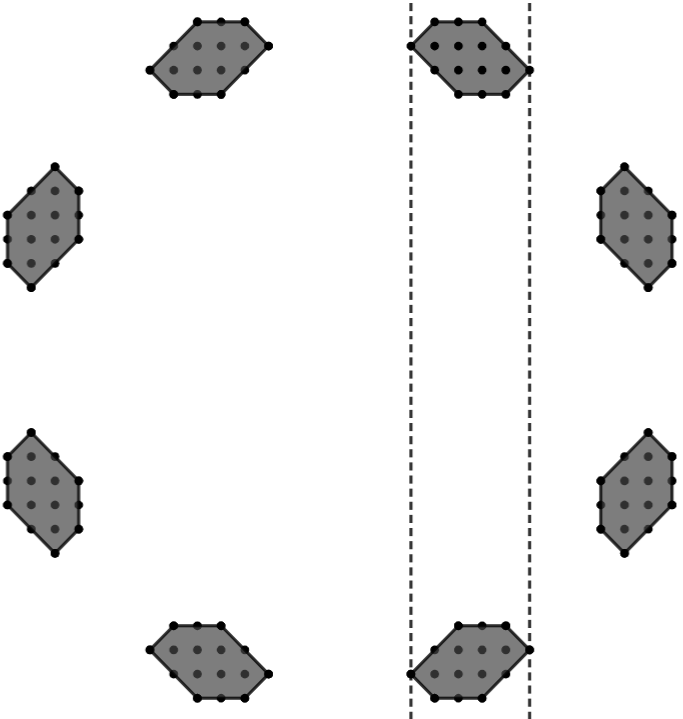}
\caption*{(b) $W(L(P_{\lambda}))$ with all the root lines containing points not in $W(L^2(P_{\lambda}))$. }
\end{subfigure}
\end{figure}

\begin{figure}[h!]

\begin{subfigure}[h!]{0.4\linewidth}
\hspace{-1cm}
\includegraphics[width = 8cm ,height = 7.8cm]{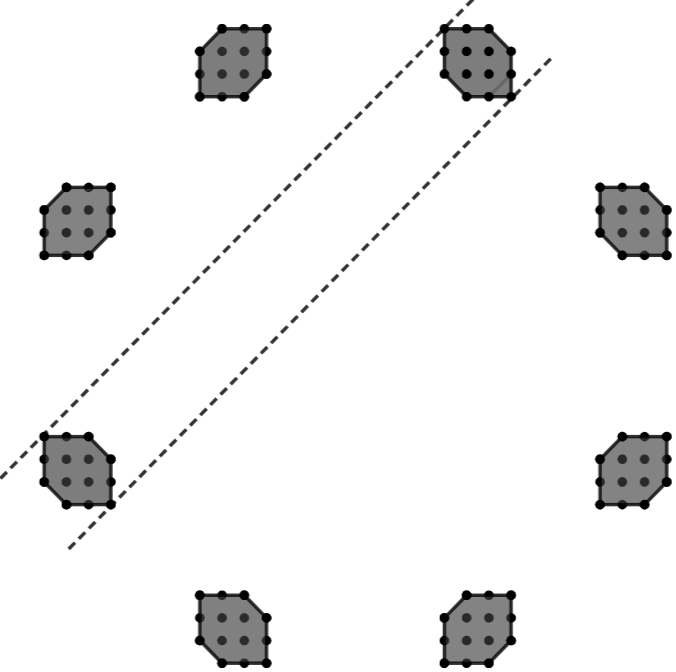} 
\caption*{(c) \ $W(L^2(P_{\lambda}))$ shown with root lines containing points not in $W(L^3(P_{\lambda}))$.}
\end{subfigure}
\hspace{2cm}
\begin{subfigure}[h!]{0.4\linewidth}
\vspace{2mm}
\includegraphics[width = 8.3cm ,height = 8.1cm]{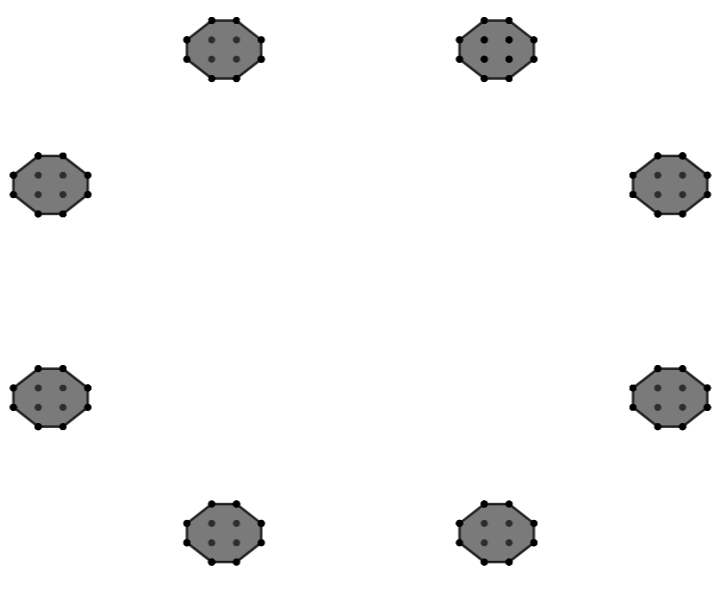}

\caption*{(d) $W(L^3(P_{\lambda}))$ and equivalently $W_{X_0}(\lambda)$ which proves Theorem \ref{GeneralizedWeylCharacter} }
\vspace{5mm}
\end{subfigure} 
\end{figure}

\newpage
Let $L$ be a line segment in $\mathbb{R}^d$ and $\Omega \subset \R^d$. We say that $\Omega$ \textit{dominates} $L$ if for any line $L'$ intersecting $\Omega$ and parallel to $L$, there exists a line segment with endpoints in $\Omega \cap L'$ whose length is greater than or equal to $|L|$. To prove Theorem \ref{GeneralizedWeylCharacter}, we first prove the following lemma regarding finite sets dominating line segments.

\begin{lemma}\label{lem:convdom}
If a finite set $S$ dominates a line segment $L$, then conv$(S)$ dominates $L$.
\end{lemma}
\begin{proof}
Assume $S$ dominates $L$ and denote the polytope conv$(S)$ by $Q$. Let $L'$ be a line parallel to $L$ which intersects $Q$. Let $L_s= L' \cap Q$. We wish to prove $|L_s| \geq |L|$. Let $x,y$ be the endpoints of $L_s$. If $x \in S$ or $y \in S$, then by domination of $S$ over $L$ the result follows.  \\

Thus, suppose $x,y \notin S$. Let $F_x$ be the minimal face of $Q$ containing $x$ and let $V$ be the vertices of $F_x$. Since each vertex $v \in V$ also lies on the boundary of $F_x$, there exists a line segment $L_v$ parallel to $L_s$ which connects $v$ to a point $v' \in Q$ and is of length $|L|$. Fix one such $\tilde{v} \in V$, and let $T_{\tilde{v}}$ be the translation map sending $\tilde{v}$ to $x$. $T_{\tilde{v}}$ sends $L_{\tilde{v}}$ to a segment of length $|L|$ connecting $x$ to some $x'$. Since $x$ is a convex combination of $V$, $x'$ is a convex combination of $\{v'| v \in V \}$ and thus $x' \in Q$. Thus, the line segment connecting $x$ and $x'$ lies in $Q$. This segment is of length $|L|$ and is a subsegment of $L_s$. 


\end{proof}

Oleg Chalykh proved the following result in \cite{C}. 

\begin{lemma}\label{lem:Oleg}
Let $L_1, \dots, L_m$ be a set of non-parallel line segments in $\mathbb{R}^d$ and $\Omega$ a convex set which dominates $L_1, \dots , L_m$. Then $\Omega$ can be represented as $\Omega = \Omega' + N$ for some convex domain $\Omega'$ where $N = L_1 + L_2+ \dots + L_m$.
\end{lemma}

We can now prove Theorem \ref{GeneralizedWeylCharacter}:

\begin{proof}[Proof of Theorem~\ref{GeneralizedWeylCharacter}] 
Suppose $m$ is the smallest integer such that $L^m(P_{\lambda}) = L^{m+1}(P_{\lambda})$. Since $P_{\lambda}$ lies in the open fundamental Weyl chamber $C$ and $W_{X_0}(\lambda)^* \subset P_{\lambda}$, we have $W_{X_0}(\lambda)^* \subset C$. It is not difficult to see that $L^m(W_{X_0}(\lambda)^*) = W_{X_0}(\lambda)^*$, so $W_{X_0}(\lambda)^* \subseteq L^m(P_{\lambda})$. We will prove the reverse inclusion. For each $r \in R^+$, we define the line segment

\begin{equation}
    \ell_{r} := \text{conv}(0, k_{r} r)
\end{equation}

 By definition of $m$, $L^m(P_{\lambda})$ dominates $\ell_{r}$ for every $r \in R^+$. By Lemma \ref{lem:convdom}, conv$(L^m(P_{\lambda}))$ also dominates each $\ell_{r}$. Lemma \ref{lem:Oleg} implies the existence of a convex set $\Omega$ such that 
\begin{equation}
    \text{conv}(L^m(P_{\lambda})) = \Omega + \sum_{r \in R^+} \ell_{r}.
\end{equation}

It suffices to prove that $\Omega = \{\lambda\}$. If $\Omega$ contains some $\mu \neq \lambda$, then      
\begin{align}
    \mu + \sum_{r \in R^+} \ell_r & \subset \text{conv}(P_{\lambda}).
\end{align}

This means $\mu + \rho$ also lies in conv$(P_{\lambda})$. We have that     
\begin{equation}
    \text{conv}(P_{\lambda}) = (\lam+{\rm cone}(R^+))\cap (\lam+\rho-{\rm cone}(R^+)). 
\end{equation}

Since $\text{conv}(P_{\lambda}) \subset \lambda + \rho - \text{cone}(R^+)$, we have $(\mu - \lambda) \in -\text{cone}(R^+)$. On the other hand, we also have that $\text{conv}(P_{\lambda}) \subset \lambda + \text{cone}(R^+)$, which implies $(\mu - \lambda) \in \text{cone}(R^+)$. Since cone$(R^+) \cap   -$cone$(R^+) = \{ 0\}$, we must have $\mu = \lambda$.
\end{proof}

\section{Weyl Spline Representations of Quasi-Interpolation Operators}
Given a finitely supported sequence $a \in S(\Lambda)$, an operator $Q_a$, defined on some subspace of continuous functions, is a \textit{local discrete quasi-interpolation operator} if (see \cite{BoxSplines}) 
\begin{equation}
    Q_a(f) := \sum_{\mu \in \Lambda}(a \ast f|_{\Lambda})(\mu)B_X(\cdot - \mu).
\end{equation}

Approximation properties of these quasi-interpolants are often stated when they are maps from an appropriate Sobolev space into the closure $\overline{\mathcal{S}(X)}$ of $\mathcal{S}(X)$ under $||\cdot||_{L^2(\R^d)}$. In this case, we have the celebrated Strang-Fix conditions (see \cite{BoxSplines}) which give us a result we state below after introducing a few required notions.  \\

For a positive integer $n$, we denote the algebra of polynomials in $n$ variables with coefficients in $\R$ as $\Pi_{n}(\R)$. If $Q_a$ is a quasi-interpolation operator on some space of functions on $\R^d$, we say that $Q_a$\textit{ reproduces polynomials up to order $m$} if for any polynomial $p \in \Pi_d(\R)$ of order at most $m$ (degree at most $m-1$) $Q_a(p)$ is well-defined and  

\begin{equation}
    Q_a(p) = p.
\end{equation}

We define $\kappa(X)$ to be the size of the smallest subset you can remove from $X$ for which the result no longer spans $\R^d$ i.e:
 \begin{equation}
 \kappa(X) := \min \{ \vert{}Y\vert{} \mid \ Y \subset X \ \text{and} \ \text{span}(X \setminus Y) \neq \mathbb{R}^d \}.
 \end{equation}

While the polynomial degree of $B_X$ is $|X|-d$, $B_X$ reproduces all polynomials of order $\kappa(X)$ and not all polynomials of order $\kappa(X)+1$. In other words, a polynomial $p$ of order $\kappa(X)+1$ exists for which there is no sequence $g: \Lambda \to \R$ such that 

\begin{equation}
    B_X \ast' g = p.
\end{equation}

\begin{theorem}\label{Deboor_Schoenberg}
    Suppose $Q$ is a local discrete quasi-interpolation operator and that $m$ is an integer such that $1 \leq m\leq \kappa(X)$. If $Q$ satisfies the following conditions: 
    \begin{itemize}
        \item $Q: H^m(\R^d) \to L^2(\R^d)$ is bounded
        \item  $Q$ reproduces all $p \in \Pi_d(\R)$ up to order $m$,
    \end{itemize}
    then $Q_h := T_hQT_h^{-1}$ satisfies, for all $f \in H^m(\R^d)$ with $C$ independent of $f$ and $h$,
    \begin{equation}
        ||Q_{h}(f) - f ||_{L^2(\R^d)} \leq Ch^m|f|_{H^m(\R^d)}.
    \end{equation} 
    \end{theorem}

If $Q$ satisfies the conditions of this theorem with $m = \kappa(X)$, then $Q$ is said to have \textit{optimal approximation order}. \\

For each $\lambda \in \Lambda^*$, we let $\mathcal{B}_{\lambda}$ be a basis of Weyl splines for $\mathcal{S}_X(\lambda)$ and let

\begin{equation}
    \mathcal{B} := \bigcup_{\lambda \in \Lambda^*} \mathcal{B}_{\lambda}.
\end{equation}
$\mathcal{B}$ is a \textit{Weyl spline basis} for $\mathcal{S}(X)$ and we denote by $\mathcal{B^*}$ its dual basis under $\langle \cdot , \cdot \rangle_X$. Let $a \in S(\Lambda)$. Since any $f \in \mathcal{S}(X)$ may be expressed in terms of the Weyl spline basis as

\begin{equation}
    f = \sum_{E \in \mathcal{B}} (-1)^{|X_0|}\langle c_f , D_XE^* \rangle_{\ell^2(\Lambda)}E 
\end{equation}
and $f \mapsto c_f$ is bijective between $\mathcal{S}(X)$ and $S(\Lambda)$, we may define $P_a: C_c(\R^d) \to \mathcal{S}(X)$ as
\begin{equation}\label{PreApproximationMap}
    P_a(f) := \sum_{E \in \mathcal{B}}  (-1)^{|X_0|}\langle a \ast f|_{\Lambda}, D_{X}E^* \rangle_{\ell^2(\Lambda)} E. 
\end{equation}

Lemma \ref{SignificantReduction} tells us that each coefficient above depends on the value of $f$ on a uniformly bounded number of points. For a generic choice of coefficient, this bound is precisely

$$|Z(X_0) + \frac{\rho}{2} \cap \Lambda| \cdot|W|$$ 
points.\footnote{Related is that the number of non-zero terms in \ref{PreApproximationMap} is bounded above by $|P(X_0)| \cdot|W|\cdot|\text{supp}(f|_{\Lambda})|$.} The next theorem tells us that $P_a(f) = Q_a(f)$ if $f \in C_c(\R^d)$ and gives the coefficients of $P_a(f)$ in \eqref{PreApproximationMap} in terms of generalized derivative samples for sufficiently smooth $f$. 

\begin{theorem}\label{Basic_Quasi_Interpolation} 
$P_a$ satisfies the following two equations:
\begin{itemize}
    \item $ \displaystyle{P_a(f) = B_X \ast' (a \ast f|_{\Lambda})}$.
    \item If $f \in C_c^{|A|}(\R^d)$, then $\displaystyle{P_a(f) = \sum_{E \in \mathcal{B}} (-1)^{|X_0|-|A|}\langle D_Af \ast' a, D_{X \backslash A  } E  \rangle_{L^2(\R^d)}E^*.}$\footnote{When $D_{X\backslash A}E$ is not necessarily in $L^2(\R^d)$, we interpret $\langle D_Af, D_{X \backslash A}E \rangle_{L^2(\R^d)}$ in a distributional sense.}
\end{itemize}
\end{theorem}

\begin{proof}
Since $\mathcal{B}$ and $\mathcal{B^*}$ are dual under $\langle \cdot , \cdot \rangle_X$, their corresponding $B_X$ and $T_X$ representations,  $\{c_E\}_{E \in \mathcal{B}}$ and $\{D_XE^*\}_{E \in \mathcal{B}}$, form, up to a sign, a biorthogonal system for $S(\Lambda)$ under $\langle \cdot , \cdot \rangle_{\ell^2(\Lambda)}$. By expressing each $E \in \mathcal{B}$ in terms of its $B_X$-representation, $B_X \ast' c_E$, we get that

\begin{align}
    P_a(f) & = B_X \ast' \sum_{E \in \mathcal{B}} (-1)^{|X_0|}\langle  f|_{\Lambda} , D_{X}E^*\rangle_{\ell^2(\Lambda)} c_{E} \\
    & = B_X \ast' (a \ast f|_{\Lambda}).
\end{align}
Now assume $f \in C_c^{|A|}(\R^d)$. We may convolve both sides of \eqref{Derivative_of_Box_Spline} by $f$, evaluate at $0$, and apply integration by parts to obtain that for any $A \subset X$,

\begin{equation}
    (\Delta_Xf)(0) = (-1)^{|A|}\int_{\R^d} D_A f D_{X \backslash A} B_X.
\end{equation}
Applying this to each coefficient of $E$ in the definition of $P$, we obtain
\begin{align}\label{Generalized_D_Samples}
    (-1)^{|X_0|}\langle a \ast f|_{\Lambda}, D_XE   \rangle_{\ell^2(\Lambda)} =
     (-1)^{|X_0|-|A|}\int_{\R^d} (D_Af \ast' a) D_{X \backslash A}E . 
    \end{align}
\end{proof}

Now, recall that a root system $R$ is \textit{irreducible} if it cannot be partitioned into two orthogonal disjoint non-empty subsets. We may now state one of the main results of this section. \\

 
\begin{theorem}\label{Section 3.2 Main Result}
    If the root system underlying $X$ is irreducible, then $a\in S(\Lambda)$ can be chosen so that $Q_a$ has optimal approximation order.   
\end{theorem}

It suffices to prove if the root system underlying $R$ is irreducible, then $a$ can be chosen so that $Q_a$, when viewed as an operator on $H^{\kappa(X)}(\R^d)$, satisfies the conditions of Theorem \ref{Deboor_Schoenberg}. \\  

We first prove that $Q_a:  H^{\kappa(X)}(\R^d) \to \overline{\mathcal{S}(X)}$ is bounded. For any $f \in C(\R^d)$, define
\begin{equation}
    \mathcal{R}(f) := f|_{\Lambda}.
\end{equation}
Thus, we represent $Q_a$ as the composition of maps 
\begin{equation}\label{composition}
    H^{\kappa(X)}(\R^d) \to_{\mathcal{R}} \ell^2(\Lambda) \to_{(B_X \ast' a) \ast' } \overline{\mathcal{S}(X)}  
\end{equation}

To see the boundedness of the second map, we recall the fact that $B_X$ forms a partition of unity, up to a constant multiple determined by $\Lambda$ (see \cite{BoxSplines}). More specifically, 

\begin{equation}
   \sum_{\lambda \in \Lambda} B_X(\cdot - \lambda) = \frac{1}{\text{det}(\Lambda)}. 
\end{equation}

We then obtain boundedness of the second map in \eqref{composition} by applying a Cauchy-Schwarz argument that yields

\begin{equation}
    ||B_X \ast' c ||_{L^2(\R^d)} \leq \frac{||c||_{\ell^2}}{|\text{det}(\Lambda)|^{1/2}}.
\end{equation}

Thus, we only need to prove the boundedness of the first map in  \eqref{composition}. It is well known that if $f \in H^s(\R^d)$ and $s > d/2$, then $\mathcal{R}$ is bounded. So, we argue that $\kappa(X) > d/2$. This is the only place we require that $X$ comes from an irreducible root system and the result follows from the following lemma.

\begin{lemma}\label{Irreducibility}
    If $R$ is an irreducible root system in $\R^d$, then $\kappa(R^+) \geq d$. 
\end{lemma}
\begin{proof}
Assume $R$ is irreducible. Let $Y \subset R^+$ such that $R^+ \backslash Y$ does not span $\R^d$ and let 
\begin{equation}
S := R^+ \backslash Y.
\end{equation}
Let $H_S:= \text{span}(S)$ and $O(H_S):= \text{span}(R \backslash H_S)$. We show that if dim$(O(H_S)) < d$, then $H_S$ and $O(H_S)\backslash H_S$, which is necessarily non-empty, are orthogonal to one another. This implies a partition of the roots into two disjoint mutually orthogonal pieces thereby contradicting the irreducibility of $R$. Since we certainly have $O(H_S) \subset \text{span}(Y)$, this would imply $\text{span}(Y) = \R^d$ giving us that  $|Y| \geq d$.     \\

So, suppose $\text{dim}(O(H_S)) < d$ and let $\alpha, \beta$ be roots where $\alpha \in R\backslash O(H_S)$ and $\beta \in R \backslash H_S$. This implies $\beta \in O(H_S)$. If $s_{\beta}(\alpha) \in O(H_S)$, then since 

\begin{equation}
s_{\beta}(\alpha) = \alpha - 2\frac{\alpha \cdot \beta}{||\beta||^2}\beta, 
\end{equation}
we may rearrange the equation above and get that $\alpha \in R \backslash O(H_S)$ is a linear combination of two elements in $O(H_S)$ which is not possible. Thus, $\alpha \cdot \beta = 0$ which implies that $R\backslash O(H_S)$ and $R\backslash H_S$ are orthogonal to each other. In particular, $R \backslash O(H_S)$ is orthogonal to $O(H_S)$ giving us the result.
\end{proof}

It follows that, for the choices of $X$ in this paper and $a \in S(\Lambda)$, $Q_a$ is bounded. To prove Theorem \ref{Section 3.2 Main Result}, we need only show that we can choose $a \in S(\Lambda)$ so that $Q_a$ reproduces polynomials up to order $\kappa(X)$. We give a proof in the appendix as the argument is standard.   \\

It may be of interest to produce an optimal quasi-interpolation operator which is also a projection onto $\mathcal{S}(X)$. Towards this, we recall that a multiset $Y \subset \R^d$ whose integer-span is a full rank lattice is \textit{unimodular} if any $d\times d$ submatrix $A \subset Y$ satisfies det$(A) \in \{-1,0,1 \}$. It is known that $X$ is unimodular if and only if the underlying root system of $X$ is the $A_n$ root system (see \cite{mcSwiggen2021}). We assume the underlying root system is $A_n$ for the rest of this section. \\

With formal power series, we define the infinite order differential \textit{Todd operator} as

\begin{equation}
    \mathcal{T}_X := \prod_{v \in X} \frac{D_v}{2\cdot\text{sinh}(D_v/2)},   \hspace{10mm}\text{where} \ \sinh(x) = \frac{e^{x}  -e^{-x}}{2} .
\end{equation}

Let $\mathcal{T}_X^*$ be the truncation of the formal power series defining $\mathcal{T}_X$ at order $|X|-d$. The following is proven about $\mathcal{S}(X)$ in \cite{Vergne} (see \cite{mcSwiggen2021} for details on the formulation stated here). 
\begin{theorem}
    Let $\eta$ be a vector in the positive cone generated by $X$ such that $\eta$ does not lie in any hyperplane spanned by elements of $X$. For any $f \in \mathcal{S}(X)$, we have
    \begin{equation}
   \widehat{\mathcal{T}_X^*}(f) := \bigg(\lim_{t \to 0}\mathcal{T}_X^*(f)(\cdot + t\eta)\bigg)\bigg|_{\Lambda} = c_f.\footnote{The choice of $\eta$ and the limit allows us view $\widehat{T_X^*}$ as the restriction to $\Lambda$ as the limit of an ordinary differential operator which does not enter the boundaries of the domains of polynomiality of $f$.}
\end{equation}
\end{theorem}

For sufficiently smooth functions, $\widehat{\mathcal{T}}_X^*(f)$ is simply a linear combination of $\{D_Af|_{\Lambda} \}_{A\subset X}$. If
\begin{equation}
    \mathcal{H} := C^{|X|-d}(\R^d) \oplus \mathcal{S}(X),
\end{equation}

we may extend the domain of $\langle \cdot, \cdot \rangle_X$ to a bilinear form on $ \mathcal{H} \times \mathcal{S}(X)$ defined as 
\begin{equation}\label{}
    (f,g)_X := (-1)^{|X_0|} \langle \widehat{\mathcal{T}_X^*}(f), d_g \rangle_{\ell^2(\Lambda)}.
\end{equation}

We will use $(\cdot,\cdot)_X$ to form a projector on $\mathcal{H}_c := \mathcal{H} \cap C_c(\R^d)$. Before doing so, recall that that if $\{D_i\}_{i \in I}$ is a finite collection of homogeneous differential operators and, for each $i \in I$,  $c_i \in S(\Lambda)$, then an operator of the form
\begin{equation}
   H(f) : = \sum_{\lambda \in \Lambda} \bigg(\sum_{i \subset I} (c_{i}\ast (D_if)|_{\Lambda})(\lambda) \bigg) B_X(\cdot - \lambda)
\end{equation}
is a \textit{local discrete differential quasi-interpolation} operator. The following map 
\begin{equation}
    \mathcal{T}(f) := \sum_{E \in \mathcal{B}} ( \widehat{\mathcal{T}_X^*}(f),E )_X E^*
\end{equation}
is well-defined, for the same reason $Q_a$ is, and is a projector on $\mathcal{H}_c$. $\mathcal{T}$ is a differential quasi-interpolation operator for the same reason $P_a$ is a local discrete quasi-interpolation operator given in  Theorem \ref{Basic_Quasi_Interpolation}. If $D$ is a differential operator and $f$ is sufficiently smooth, we have  
\begin{equation}
    \langle Df, D_XE \rangle_{\ell^2(\Lambda)} = (-1)^{|A|}\langle D_ADf , D_{X \backslash A} E \rangle_{L^2(\R^d)},  
\end{equation}
so the coefficients of $\mathcal{T}(f)$ are expressible as generalized derivative samples of $f$ induced by $\mathcal{B}$. To prove $\mathcal{T}$ has optimal approximation order, we define $\mathcal{T}^*$ on $C^{|X|-d}(\R^d)$ as 
\begin{equation}
    \mathcal{T}^*(f) := B_X \ast' \widehat{\mathcal{T}_X^*(f)}. 
\end{equation}
$\mathcal{T}^*(f) = \mathcal{T}(f)$ for $f \in \mathcal{H}_c$ for the same reason $P_a(f) = Q_a(f)$ for $f \in C_c(\R^d)$.

\begin{theorem}\label{DQIO}
 Let $m \geq |X|-d +\kappa(X)$. $\mathcal{T^*}$ has optimal approximation order on $H^m(\R^d)$ in the sense that if $\mathcal{T}_h := T_h\mathcal{T}T_h^{-1}$ and $f \in H^m(\R^d)$, then there exists $C > 0$ such that
   \begin{equation}
       ||\mathcal{T}_h(f) - f ||_{L^2(\R^d)} \leq Ch^{\kappa(X)}|f|_{H^{\kappa(X)}(\R^d)}.
   \end{equation} 
\end{theorem}
\begin{proof}
      Since $\mathcal{T}$ is the identity on $\mathcal{S}(X)$, $\mathcal{T}$ locally reproduces all $p \in \Pi_d$ up to order $\kappa(X)$. Thus, $\mathcal{T}^*$ reproduces all polynomials in $\Pi_d$ up to order $\kappa(X)$. Since $D_A$ is bounded from $H^m(\R^d)$ to $H^{\kappa(X)}(\R^d)$ for all $A \subset X$ and $\mathcal{T}_X^*(f)$ is a linear combination of $\{D_Af\}_{A \subset X}$, $\mathcal{T}_X^*$ is bounded from $H^m(\R^d)$ into $H^{\kappa(X)}(\R^d)$. \\
     
      Thus, $\mathcal{T}^*$ is bounded from $H^m(\R^d)$ into $L^2(\R^d)$ since $\mathcal{R}(f)$ is bounded from $H^s(\R^d)$ into $\ell^2(\Lambda)$ for any $s > d/2$ and $B_X \ast'$ is bounded from $\ell^2(\Lambda)$ into $L^2(\R^d)$. Let $\Pi_{d,\kappa(X)}$ be all $p \in \Pi_d(\R)$ of order $\kappa(X)$ or less. Let $K$ be a compact set and $K^*$ be an enlargement of $K$ containing all points needed to compute $\mathcal{T}(f)|_{K}$. Polynomial reproduction and boundedness of $\mathcal{T}$ imply that
    \begin{align}
        ||\mathcal{T}(f) - f||_{L^2(K)} & \leq \inf_{p \in \Pi_{d,\kappa(X)}} ||\mathcal{T}(f -p)||_{L^2(K)} + ||p- f||_{L^2(K)} \\
        & \leq Ch^{\kappa(X)}\inf_{p \in \Pi_{d,\kappa(X)}}||f - p||_{H^m(K^*)} \\
        & \leq Ch^{\kappa(X)}|f|_{H^{\kappa(X)}(K^*)}
    \end{align}
     where the last inequality comes from the Bramble-Hilbert lemma. Due to translation invariance of $\Pi_{d,\kappa(X)}$, the constant $C$ is uniform over translations of $K^*$. The rest of the proof follows from standard arguments using the dilation invariance of $\Pi_{d,\kappa(X)}$ and summing over translations of an appropriately chosen $K$ and $K^*$. 
\end{proof}

Other approximation maps can be constructed from $\langle \cdot, \cdot \rangle_X$ either by extending the domain $\langle \cdot, \cdot\rangle_X$ directly and building a projector or from forming the analysis-synthesis operator coming from an inner product similar to $\langle \cdot, \cdot\rangle_X$. It is worth emphasizing that the especially novel part of this construction is that the coefficients of $P_a$ and $\mathcal{T}$, in their Weyl spline representations, can be expressed as generalized derivative samples.


\section{Constructing Biorthogonal Weyl Spline Systems}

In this section, we look at how to compute $\langle \cdot , \cdot \rangle_X$ on $\cS_X(\lambda)$ and subsequently how to construct a biorthogonal system of Weyl splines. We begin with the following lemma which is critical for understanding $\cS_X(\lambda)$ and ultimately tells us how to construct a Weyl spline basis. 

\begin{lemma}\label{Weyl_Spline_Index}
For any $\lambda \in \Lambda^*$, the map $E \mapsto c_E|_{W(\lambda)}$ is a linear isomorphism between $\mathcal{S}_X(\lambda)$ and the space of functions $\{f| f:W(\lambda) \to \R \}$. 
\end{lemma}

\begin{proof}

We first show dim$(\cS_X(\lambda)) = |W(\lambda)|$. Define the space $V$ as 
\begin{equation}
V := \text{span}\{B_X(\cdot -\mu)\}_{\mu^* \leq \lambda}. 
\end{equation}
For each $\mu \in \Lambda$ consider the map on $V$:
\begin{equation}
    H_{\mu}: E \mapsto (-1)^{|X_0|}(D_XE)(\mu).
\end{equation}

$H_{\mu}$ is a linear functional on $V$ and Property \ref{property 2} in the definition of $\mathcal{S}_X(\lambda)$ implies that  

\begin{equation}\label{SomeKernels}
\mathcal{S}_X(\lambda) = \bigcap_{\{\mu \in \Lambda| \ \mu^* < \lambda\}} \text{ker}( H_{\mu} ). 
\end{equation}

The set of linear functionals $\{H_{\mu}\}_{\mu^* \leq \lambda}$ is linearly independent since $\{B_X(\cdot - \mu) | \mu^* \leq \lambda \}$ is linearly independent and we have the representation
\begin{equation}
    H_{\mu}(E) = \langle B_X(\cdot - \mu), E\rangle_X.
\end{equation} 

Thus, \eqref{SomeKernels} and basic linear algebra imply dim$(\mathcal{S}_X(\lambda)) = |W(\lambda)|$\footnote{This statement follows trivially from Opdam's orthogonal basis construction; however, we presented a self-contained direct proof of this.}. 

We now show that the map $E \mapsto c_E|_{W(\lambda)}$ is injective. If $E \in \mathcal{S}_X(\lambda)$ and $c_E|_{W(\lambda)} = 0$, then 
\begin{equation}
    \text{supp}(c_E) \subset \{\mu \in \Lambda| \ \mu^* < \lambda \}.
\end{equation}
By the definition of $S_X(\lambda)$, we have $D_XE = 0$ on the set above. Consequently, 
\begin{equation}
    \langle E, E \rangle_X = (-1)^{|X_0|}\langle c_E,d_E \rangle_{\ell^2(\Lambda)}  = 0.
\end{equation}

Since $\langle \cdot, \cdot \rangle_X$ is an inner product, this implies $E = 0$.
\end{proof}

This lemma implies that we may completely determine $E \in \mathcal{S}_X(\lambda)$ by the value of $c_E$ on $W(\lambda)$. By Lemma \ref{Weyl_Spline_Index}, we are permitted to associate each $f: W(\lambda) \to \R$, with the unique spline $E_f \in \mathcal{S}_X(\lambda)$ that satisfies   
\begin{equation}
c_{E_f}|_{W(\lambda)} = f.\footnote{For a generic choice of $\lambda \in \Lambda^*$, $W(\lambda)$ can be identified with $W$; however, if $\lambda$ lies on the boundary of the fundamental Weyl chamber, $|W(\lambda)| < |W|$. } 
\end{equation}

From here on, we assume $\lambda$ satisfies the conditions of Theorem \ref{GeneralizedWeylCharacter}. For each such $\lambda$, we refer to the vertices of $Z_{X_0}(\lambda)$ as the \textit{extreme spectral points} of $W_{X_0}(\lambda)$ which we denote as 

$$V(Z_{X_0}(\lambda)).$$ 

$V(Z_{X_0}(\lambda))$ can be understood in terms of the Weyl group. To see this, observe that the vertices of $Z(X_0)$ may be indexed by $W$ since the vertex set of $Z(X_0)$ is  

\begin{equation}
\{\tau(\rho/2)| \ \tau \in W \}.
\end{equation}

$Z_{X_0}(\lambda)$ is the union of $|W|$ disjoint zonotopes with each zonotope lying in an orbit of the fundamental Weyl chamber $C$. We may index each such zonotope by a unique $\sigma \in W$. Therefore, we may index each point in $V(Z_{X_0}(\lambda))$ by a unique $(\sigma , \tau) \in W \times W$ which corresponds to the translation of $\tau(\rho/2)$ lying in the chamber $\sigma(C)$. We write this precisely with the following notation:

\begin{equation}\label{SpectralIndexing}
 (\sigma, \tau)_{\lambda} := \sigma\bigg( \lambda + \frac{\rho}{2} + \tau\bigg(\frac{\rho}{2}\bigg) \bigg). 
 \end{equation}

We say $(\sigma, \tau)_{\lambda}$ and $(\sigma, w)_{\lambda}$ are \textit{adjacent} if the line segment connecting them is an edge of 
\begin{equation}
    \sigma \bigg(\lambda + \frac{\rho}{2} + Z(X_0)\bigg).
\end{equation}

If $w \in W$ satisfies that $(\sigma , \tau)_{\lambda}$ and $(\sigma , w \tau)_{\lambda}$ are adjacent, we say that $w$ is \textit{permissible} at $(\sigma , \tau)_{\lambda}$. It's easy to see that if $w$ is permissible at any $(\sigma , \tau)_{\lambda}$, $w$ is necessarily a reflection. It follows from $\eqref{SpectralIndexing}$ that if $w$ is permissible at $(e,\tau)_{\lambda}$, then 

\begin{equation}
(e,\tau)_{\lambda}, \ (e, w\tau)_{\lambda}, \ (w, w \tau)_{\lambda}, \ (w, \tau)_{\lambda}
\end{equation}
are collinear. By applying $\sigma \in W$ to the line passing through these points, we obtain that 

\begin{equation}\label{colinear}
(\sigma, \tau)_{\lambda}, (\sigma, w \tau)_{\lambda}, (\sigma w , w \tau)_{\lambda} , (\sigma w, \tau)_{\lambda}
\end{equation}
are also collinear\footnote{It's important to note that $w$ may not reflect $(\sigma, \tau)_{\lambda}$ in the same direction that $w$ reflects 
$(e,\tau)_{\lambda} $.}. We will see later that the collinearity in \eqref{colinear} and properties of truncated powers allow us to express the value of $D_XE$ on $(\sigma, w \tau)_{\lambda}$ and $(\sigma w , w \tau)_{\lambda}$ in terms of $(\sigma, \tau)_{\lambda}$ and $(\sigma w, \tau)_{\lambda}$. Thus, if we have knowledge of $D_XE$ on $\{(\sigma, \tau)_{\lambda} | \sigma \in W \} $ for some $\tau \in W$, we can find the value of $D_XE$ on adjacent points of this set. \\

Recall the \textit{simple roots}, $\Pi \subset R^+$, defined in Appendix $A.2$. In $Z(X_0)$, the points adjacent to $\rho/2$ are $\{s_{\alpha}(\rho/2)| \alpha \in \Pi \}$. Thus, the vertices adjacent to $\tau(\rho/2)$ are $\{\tau s_{\alpha}(\rho/2)| \ \alpha \in \Pi  \}$. By translating and applying $\sigma \in W$, as in \eqref{SpectralIndexing}, we see the adjacent points to any $(\sigma,\tau)_{\lambda}$ are $\{ (\sigma,\tau s_{\alpha})_{\lambda} | \alpha \in \Pi \}$. For any $\tau \in W$ and $\alpha \in R$, the identity

\begin{equation}
\tau s_{\alpha} \tau^{-1} = s_{\tau(\alpha)}
\end{equation}

is well known. So, we may equivalently represent the points adjacent to $(\sigma , \tau)_{\lambda}$ as 

\begin{equation}\label{adjacent points}
\{(\sigma , s_{\tau(\alpha)} \tau )_{\lambda}| \  \alpha \in \Pi  \}.
\end{equation}
It's not difficult to see that, for any $\alpha \in \Pi$, $s_{\alpha}$ permutes $R^+\backslash \{\alpha \}$ and sends $\alpha $ to $-\alpha$. Since $\tau(R^+)$ is a choice of positive roots for $R$ with $\tau(\Pi)$ being its corresponding set of simple roots, we have that for any $\alpha \in \Pi$

\begin{align}
    s_{\tau(\alpha)}\tau(\rho/2) & = s_{\tau(\alpha)} \bigg(\frac{1}{2}\sum_{\beta \in \tau(R^+)} k_{\beta}\beta \bigg) \\
    & = \bigg(\frac{1}{2}\sum_{\beta \in \tau(R^+)} k_{\beta}\beta \bigg) - k_{\tau(\alpha)} \tau(\alpha) \\
    & = \tau(\rho/2) - k_{\alpha}\tau(\alpha).
\end{align}

\begin{lemma}\label{ImportantLemma}
    
Suppose $(r_1, \dots , r_{N})$ is a sequence of distinct positive roots which satisfies:

\begin{enumerate}
    \item $r_1 \in \Pi$
    \item $r_{i+1} \in (s_{r_{i}}\dots s_{r_1}\Pi) \backslash \{-r_1, \dots , -r_{i}  \}$ for any $1 \leq i < N$.
\end{enumerate}
If $N<|R^+|$, there exists a positive root $r_{N+1}$ for which $(r_1, \dots , r_{N+1})$ is a sequence of distinct roots also satisfying the conditions above. Finally, given such a sequence, suppose $w_0 = e$, $w_i := s_{r_i} \dots s_{r_1}$, and $k_i = k_{r_i}$ for all $1 \leq i \leq |R^+|$. For all $1 \leq i <|R^+|$,
\begin{equation}\label{lemma_equation}
w_{i+1}(\rho/2) = w_i(\rho/2) - k_{i+1}r_{i+1}. 
\end{equation}

\end{lemma}

\begin{proof}
Let $(r_1, \dots , r_N)$ be as stated in the lemma. Since $r_{j+1} \in w_j(\Pi)$ for all $j < N$, we have 

\begin{equation}\label{eq}
w_{j+1}(\rho/2) = w_j(\rho/2) - k_{j+1}r_{j+1}.
\end{equation}
Thus, \eqref{lemma_equation} holds for all $j < N$. We now show if $r \in R^+$ lies in $w_N(\Pi) \backslash\{-r_1, \dots , -r_N \}$, then $r$ is distinct from $\{r_1, \dots , r_N\}$. Suppose $r_j \in w_N(\Pi) \backslash \{-r_1, \dots , -r_N \}$ for some $1 \leq j \leq N$. This implies $w_N(\rho/2) - k_jr_j$ is a vertex of $Z(X_0)$. Since $Z(X_0)$ contains all linear combinations of vectors in $X_0$ with coefficient $1/2$ or $-1/2$ and
\begin{equation}
w_N(\rho/2) = \frac{\rho}{2} - \sum_{i=1}^{N} k_ir_i,    
\end{equation}
$w_N(\rho/2) + k_jr_j \in Z(X_0)$. However, if $w_N(\rho/2) - k_jr_j \in Z(X_0)$ and $w_N(\rho/2) + k_jr_j \in Z(X_0)$, then $w_N(\rho/2)$ cannot be a vertex of $Z(X_0)$ giving us a contradiction. \\

Now suppose $w_N(\Pi) \backslash \{-r_1, \dots , -r_N \} $ contains no positive roots. Then all roots in $w_N(\Pi)$ are negative and so,  every point adjacent to $w_N(\rho/2)$ may be expressed as 

\begin{equation}
    w_N(\rho/2) - k_{r}w_N(r)
\end{equation}
for some $r \in \Pi$. Since $Z(X_0) \subset \rho/2 +\text{cone}(-\Pi)$, we may apply $w_N$ to obtain
\begin{align}
 Z(X_0) & \subset w_N(\rho/2) + \text{cone}(w_N(-\Pi))\\ 
 & \subset w_N(\rho/2) + \text{cone}(r_1, \dots , r_{N}).
\end{align}

Thus, every point adjacent to $w_N(\rho/2)$ is obtained by adding positive roots to it. This could only be the case if $w_N(\rho/2)$ is the lowest point of $Z(X_0)$ in the dominance ordering or equivalently $w_N(\rho/2) = -\rho/2$. This occurs precisely when $N = |R^+|$.

\end{proof}

We refer to the sets $W(\lambda+ \rho)$ and $W(\lambda)$ as the \textit{outer} and \textit{inner} \textit{spectral points of $W_{X_0}(\lambda)$} respectively. From here on, we fix a sequence of positive roots  $(r_1, \dots , r_{|R^+|})$ satisfying the properties of Lemma \ref{ImportantLemma}. Using the notation of Lemma \ref{ImportantLemma}, for any $\sigma \in W$, the sequence

\begin{equation}\label{walk}
(\sigma,w_0)_{\lambda}, (\sigma , w_1)_{\lambda}, \dots , (\sigma , w_{|R^+|})_{\lambda}
\end{equation}
begins with an outer spectral point, ends with an inner spectral point, and any two consecutive points are adjacent. In lieu of this, we define for any $1 \leq i \leq |R^+|$ 
\begin{equation}
s_i := s_{r_i}. 
\end{equation}

and refer to $(s_1 , \dots , s_{|R^+|})$ as an \textit{OI sequence}.\footnote{Lemma \ref{ImportantLemma} implies that $w_{|R^+|}$ is the longest element of a Weyl group as it sends $R^+$ to $-R^+$. This element is typically denoted as $w_0$, but we denote this representation of it as $w_{|R^+|}$.} We will give a recursive formula for the value of $D_XE$ on points of the form \eqref{walk}. Since any vertex of $Z(X_0)$ lies in the set   

$$\{\rho/2, w_1(\rho/2) , \dots , w_{|R^+|}(\rho/2) \}$$
for some OI sequence, this recursive formula gives us $D_XE$ on all of $V(Z_{X_0}(\lambda))$. We need a basic property of truncated powers which we state here for completeness. 

\begin{lemma}\label{Truncated_Power_Lemma}
Let $v \in \R^d$ be a non-zero vector, $k$ be a positive integer, and $Y$ be the multi-set 
\begin{equation}
Y := [v, \dots , v] \in \R^{d \times k}.
\end{equation}

Suppose $\{t_1, \dots , t_N\}$ is a distinct set of real numbers and $\lambda_i := vt_i +b$ for $1 \leq i \leq N$ is a set of collinear points in $\R^d$. There exists coefficients $c_1 , \dots , c_N$ such that the distribution
\begin{equation}
    \sum_{i=1}^N c_iT_Y(\cdot - \lambda_i)
\end{equation}
is compactly supported if and only if $N \geq k+1$. Moreover, if $N = k+1$, then ,up to scaling,

\begin{equation}
    c_i = \prod_{j=1, j\neq i}^{k+1} \frac{1}{(t_i - t_j)}.
\end{equation}
\end{lemma}  

The result above is basic in the uni-variate setting and the proof is essentially the same so we omit it here. Using this, we can now prove the following result.

\begin{theorem}\label{MainAlgorithm}

Let $(s_1, \dots , s_{|R^+|})$ be an OI sequence for $W_{X_0}(\lambda)$ and $(r_1, \dots , r_{|R^+|})$ be the corresponding sequence of positive roots. Furthermore, define

\begin{equation}
\gamma(x,k):=(-1)^k\frac{\Gamma(x+k)\Gamma(x-k)}{\Gamma(x)^2}.
\end{equation}

For each $1\leq i \leq |R^+|$, let $k_i$ be the multiplicity of $r_i \in X_0$ and 

\begin{equation}\label{StationaryDistances}
x_i := \frac{||(\lambda + \rho/2) - s_i(\lambda + \rho/2) ||_2}{||r_i ||_2}. 
\end{equation}

Let $\gamma_i := \gamma(x_i, k_i)$. Finally suppose $E \in \mathcal{S}_{X}(\lambda)$. For any $\sigma \in W$ and $1 \leq i \leq |R^+|$:
\begin{align}
    D_XE((\sigma, w_i)_{\lambda}) & =  \gamma_i  \bigg( D_XE((\sigma, w_{i-1})_{\lambda}) -  \frac{k_i }{x_i} D_XE(( \sigma s_i , w_{i-1})_{\lambda})     \bigg).
\end{align}
\end{theorem} 

\begin{proof}

We fix $i$ where $1 \leq i \leq |R^+|$ and let $L_i$ be the line containing $(e , w_{i-1})_{\lambda}$ and $(s_{i} , w_{i-1})_{\lambda}$. We label all of the $2k_i+2$ points in $W_{X_0}(\lambda) \cap L_i$ as follows:

\begin{equation}\label{L_ilabelling}
\lambda_j := \begin{cases}
    (e, w_{i-1})_{\lambda} - (j-1)r_i \  & \text{if} \  1 \leq j \leq k_i+1 \\
     (s_i,w_{i})_{\lambda} - ((j-2)- k_i)r_i  \ & \text{if} \ k_i+2 \leq j \leq 2k_i+2.
\end{cases} 
\end{equation}

Under this labeling, Lemma \eqref{ImportantLemma} implies $\lambda_{k_i+1} = (e,w_i)_{\lambda}$. The proof we give here will express $D_XE((e,w_i)_{\lambda})$ in terms of the values of $D_XE$ on $(e, w_{i-1})_{\lambda}$ and  $(s_{i}, w_{i-1})_{\lambda}$.

 \begin{figure}[h!]
     \caption{The figure below shows part of the Weyl spectrum in the case where $R= A_2$ and $R^+$ is the same as in Figure \ref{Fig2} ordered as $r_1 = (1,0), r_2 = (\frac{1}{2}, \frac{\sqrt{3}}{2}), r_3 = (-\frac{1}{2}, \frac{\sqrt{3}}{2})$, and $X_0 = [R^+|R^+]$. We show $L_2$ from the proof which contains $(e,w_1)_{\lambda} = \lambda_1$ and $(s_2,w_1)_{\lambda} = \lambda_6$.}
     \label{fig:placeholder}
     \centering
     \includegraphics[width = 7.8cm ,height = 8.0cm]{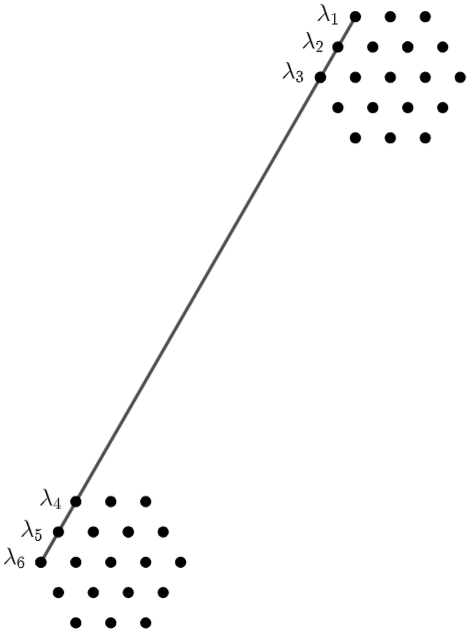}
 \end{figure}

For any $r \in X$, let $X_{r} \subset X$ be the multi-set gotten by removing all vectors from $X$ except for $r$. Define $F_i := D_{X \backslash X_{r_i}} E$. By Theorem \ref{GeneralizedWeylCharacter}, we have a representation of $F_i$ as
\begin{equation}
    F_i = \sum_{\mu \in W_{X_0}(\lambda)} d_E(\mu)T_{X\backslash X_{r_i}}(\cdot - \mu )
\end{equation}
We define the restriction of the sum above over indices lying $L_i \cap W_{X_0}(\lambda)$ as
\begin{equation}
    F_{L_i} := \sum_{j= 1}^{2k_i+2} a_jT_{X_{r_i}}(\cdot - \lambda_j).
\end{equation}

$F_{L_i}$ is compactly supported and so we may apply Lemma \ref{Truncated_Power_Lemma} to solve for the coefficients above. For each $1 \leq i \leq |R^+|$, define

\begin{equation}\label{distances}
    \overline{x_i} :=  \frac{|| (e, w_{i-1} )_{\lambda} - (s_i, w_{i})_{\lambda} ||_2}{||r_i||_2} 
\end{equation}
or equivalently, the number of lattice points along the line segment connecting $\lambda_1 = (e,w_i)_{\lambda}$ and $\lambda_{k_i + 2} = (s_i,w_i)$. If we apply Lemma \ref{Truncated_Power_Lemma} with $T_{X_{r_i}}$ and the collinear points $\{\lambda_1 , \dots , \lambda_{2k_i+1} \}$ with $c_1 = (D_XE)((e,w_{i-1})_{\lambda}))$, the resulting coefficient of $T(\cdot - \lambda_{k_i+1})$ is  

\begin{equation}
    c_{k_i+1} = D_XE((e,w_{i-1})_{\lambda}) \gamma(\overline{x_i},k_i) 
\end{equation}
If we apply Lemma \ref{Truncated_Power_Lemma} with $T_{X_{r_i}}$ and $\{\lambda_2 , \dots , \lambda_{2k_i+2} \}$ with $c_{2k_i+2} = D_XE((s_i,w_{i-1})_{\lambda})$, then

\begin{equation}
     c_{k_i+1} = -\gamma(\overline{x_i},k_i) \frac{k_i}{\overline{x_i}}D_XE((s_i,w_{i-1})_{\lambda}) 
\end{equation}

Summing the two resulting coefficients of $T_{X_{r_i}}(\cdot - \lambda_{k_i+1})$, we obtain that 


\begin{equation}
 a_{k_i+1} = \gamma(\overline{x_i},k_i)\bigg( D_XE((e, w_{i-1})_{\lambda}) -  \frac{k_i }{\overline{x_i}} D_XE(( s_i , w_{i-1})_{\lambda})     \bigg).
\end{equation}

This implies a recursive formula for $D_XE((e,w_i)_{\lambda}) = a_{k_i+1}$ for $1 \leq i \leq |R^+|$. We may reproduce the steps above with the points $\{\sigma(\lambda_1), \dots \sigma(\lambda_{2k_i+2}) \}$ to compute the coefficient of $(\sigma, w_i)$ for each $\sigma \in W$. This gives the result of the theorem with $\overline{x_i}$ in place of $x_i$. \\

We now prove that $x_i = \overline{x_i}$. We do this by using the translation invariance of the distance function. We translate both $(e,w_{i-1})_{\lambda}$ and $(s_i,w_i)_{\lambda}$ in \eqref{distances} by $-w_{i-1}(\rho/2)$. By \eqref{SpectralIndexing}, 
\begin{equation}
(e,w_{i-1})_{\lambda} - w_{i-1}(\rho/2) = \lambda + \rho/2
\end{equation}

and since $s_iw_i = w_{i-1}$, we also have  

\begin{align}
    (s_i,w_{i})_{\lambda} - w_{i-1}(\rho/2) & = s_i\bigg( \lambda + \frac{\rho}{2} + w_i(\rho/2) \bigg) - w_{i-1}(\rho/2)\\
    & = s_i(\lambda + \rho/2).    
\end{align}

This gives the equivalence of \eqref{StationaryDistances} and \eqref{distances}. 
\end{proof}

We may define an analogue of \eqref{distances} for any $r \in R^+$ and observe that 

\begin{equation}
    x_r := \frac{||(\lambda + \rho/2) - s_{r}(\lambda + \rho/2) ||_2}{||r||_2} 
    = \frac{ 2\langle \lambda + \rho/2 , r \rangle}{||r||_2^2}.
\end{equation}

This can be used to compute $x_r$ for all $r \in R^+$ in terms of $\{x_{\alpha} |\alpha \in \Pi\}$.

\begin{figure}[h!]
\caption{We consider $W_{X_0}(\lambda)$ under the same assumptions of Figure (\ref{Fig2}). We order $R^+$ as $r_1 = (1,0), r_2 = (\frac{1}{2}, \frac{\sqrt{3}}{2}) \ \text{and}, r_3 = (-\frac{1}{2}, \frac{\sqrt{3}}{2})$. The top left shows $W_{X_0}(\lambda)$ with $x_1$ and $x_2$ labeled and $x_3$ labeled as $x_1 + x_2$. Each subsequent figure shows $W_{X_0}(\lambda)$ with the root lines $\{w(L_i)| w\in W\}$ used in computing the coefficients of any point $(\sigma,w_i)_{\lambda}$. The table below lists all the coefficients computed of $D_XE_f$, in counter clockwise order, with $f = \delta_e$.   }
    \begin{subfigure}{0.2\linewidth}
    \includegraphics[width = 7.8cm ,height = 6.5cm]{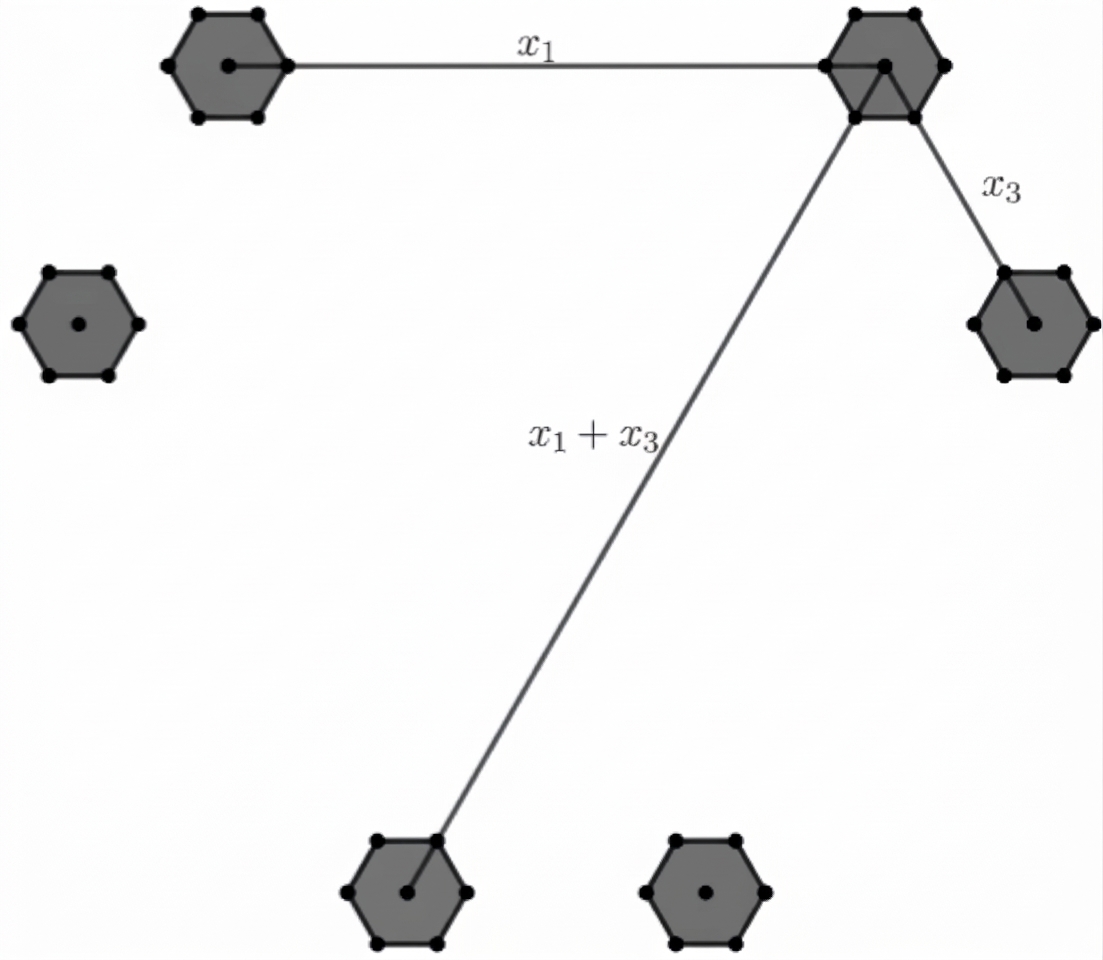} 
    \end{subfigure}
    \hspace{6cm}
    \begin{subfigure}{0.2\linewidth}
     \includegraphics[width = 7.8cm ,height = 6.5cm]{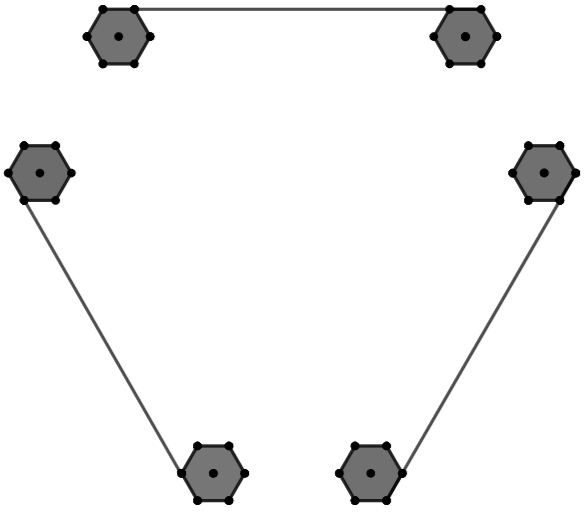} 
    \end{subfigure}
\end{figure}

\begin{figure}[h!]
    \begin{subfigure}{0.2\linewidth}
 \includegraphics[width = 7.8cm ,height = 6cm]{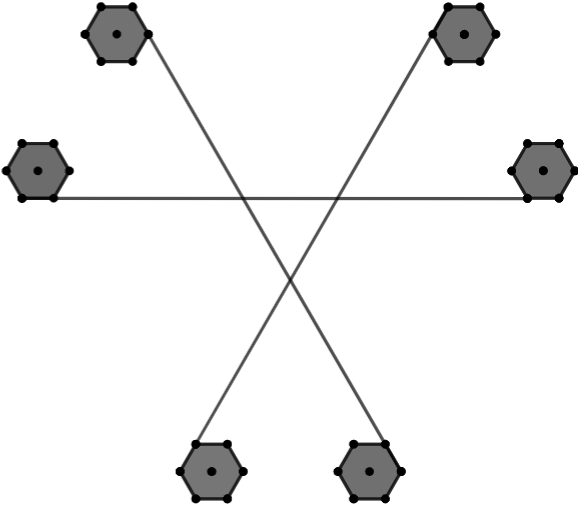}
    \end{subfigure}
    \hspace{8.5cm}
    \begin{subfigure}{0.2\linewidth}
    \includegraphics[width = 7.8cm ,height = 6cm, center]{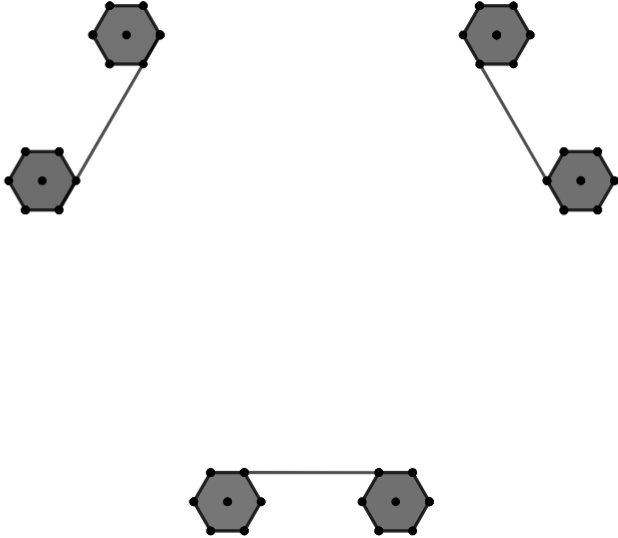}
    \end{subfigure}
\end{figure}

\begin{center}
 \begin{tabular}{|p{1.6cm}       | p{1.5cm}    |   p{4.0cm}      | p{6.5cm}    | } 
 \hline 
  $f$ &  $(W, w_1)$ &  $(W,w_2)$ & $(W,w_3)$\\ [1ex] 
 \hline 
  $1$ & $\displaystyle{ \frac{-x_1}{(x_1-1)} }$  & $\displaystyle{ \frac{x_1}{(x_1-1)} \frac{(x_1+x_3)}{(x_1+x_3-1)} }$ & $\displaystyle{ -\frac{x_1}{(x_1-1)} \frac{(x_1+x_3)}{(x_1+x_3-1)} \frac{x_3}{(x_3-1)}}$\\ [1.5ex] 
 
 \hline
 $0$ & $\displaystyle{ \frac{1}{(x_1-1)}}$ & $\displaystyle{- \frac{1}{(x_1-1)}}$$\displaystyle{ \frac{(x_1+x_3)}{(x_1+x_3-1)}}$ & $\displaystyle{ \frac{1}{(x_1-1)}}$   $\displaystyle{ \frac{(x_1+x_3)}{(x_1+x_3-1)}}$ $\displaystyle{\frac{x_3}{(x_3-1)}}$ \\ [1.5ex] 
 \hline
 $0$  &  $0$ & $0$ & $\displaystyle{ - \frac{1}{(x_1 -1)} \frac{x_1+x_3}{(x_1 + x_3 -1)} \frac{1}{(x_3-1)} } $ \\[1.5ex]
  
  \hline
  $0$ & $0$ & $\displaystyle{ \frac{x_1}{(x_1-1)}\frac{-1}{(x_1 + x_3 -1)}}$ & $\displaystyle{ \frac{(x_1x_3 +1)}{(x_1-1)(x_1+x_3-1)(x_3-1)} }$ \\[1ex]
  \hline
  $0$ & $0$ & $\displaystyle{ \frac{1}{(x_1-1)}}$$\displaystyle{ \frac{1}{(x_1+x_3-1)}}$ & $\displaystyle{- \frac{1}{(x_1-1)} \frac{x_1+x_3}{(x_1+x_3-1)} \frac{1}{(x_3-1)} }$ \\[1ex]
 \hline
 $0$ & $ \displaystyle{ 0}$ & $\displaystyle{ 0}$ & $\displaystyle{ \frac{x_1}{(x_1-1)} \frac{(x_1+x_3)}{(x_1+x_3-1)} \frac{1}{(x_3-1)} }$ \\ [1ex]
 \hline
\end{tabular}
\end{center}
\vspace{2mm}

\vspace{2mm}

From here on, we fix an OI sequence $(s_1, \dots , s_{|R^+|} )$ for $W_{X_0}(\lambda)$. We denote by $\ell^2(W)$ the space of all real-valued functions on $W$ under the Euclidean inner product $\langle \cdot, \cdot \rangle_{\ell^2(W)}$ and for each $1 \leq i \leq |R^+|$, define the operator

\begin{equation}
P_i: f \mapsto f( \cdot s_i) 
\end{equation}

for any $f \in \ell^2(W)$. Furthermore, with the notation of Theorem \ref{MainAlgorithm}, we define
\begin{equation}\label{Coefficients}
    a_i := \gamma(x_i,k_i ) \ \ \text{and} \ \ b_i := \frac{k_i}{x_i} \gamma(x_i,k_i).
\end{equation}

For any $f: W \to \R$ and $1 \leq n \leq |R^+| $, we define the map $M_n:=M_n(X,\lambda)$ on $\ell^2(W)$ via

\begin{equation}
    (M_nf)(w) := (D_X E_f)( (w,w_n)_{\lambda} ).
\end{equation} 

We interpret $M_n$ as mapping a choice of outer spectral coefficients, $f$, to the spectral coefficients of $E$ at the points $\{ (\sigma, w_n) | \sigma \in W\}$. Theorem \ref{MainAlgorithm} implies that   

\begin{equation}\label{matrixrep}
M_n = (a_{n}I - b_{n}P_{n}) \dots (a_1 I- b_1P_{1} ).
\end{equation}

This representation of $M_n$ makes it easy to retrieve the outer spectral coefficients from the coefficients of any $W$ orbit of extreme spectral points since
\begin{equation}\label{Inverse}
    M_n^{-1} = \frac{(a_1I + b_1P_{1})}{a_1^2-b_1^2} \dots \frac{(a_nI + b_nP_{n})}{a_n^2-b_n^2}.
    \end{equation}
We give special attention to $M_{|R^+|}$ and define  

\begin{equation}
M := (-1)^{|X_0|}M_{|R^+|}.
\end{equation}

$M_{|R^+|}$ sends a choice of outer spectral coefficients to its corresponding inner spectral coefficients. Theorem \ref{GeneralizedWeylCharacter} implies that for any $E_f , E_g \in \cS_X(\lambda)$, the supports of $c_{E_f}$ and $d_{E_g}$ intersect only at $W(\lambda)$. By Theorem \ref{MainAlgorithm}, for any $E_f,E_g\in \calS_X(\lambda)$, we have that 

\begin{equation}\label{spd of T}
    \langle E_f ,  E_g\rangle_{X} = \langle f, Mg \rangle_{\ell^2(W)}.
\end{equation}

Suppose $B$ is a basis for $\ell^2(W)$ and $B^*$ its dual basis under $\langle \cdot , \cdot \rangle_{\ell^2(W)}$. We  denote by $\mathcal{B}_{\lambda}$ the corresponding spline basis for $\mathcal{S}_X(\lambda)$:

\begin{equation}
\mathcal{B}_{\lambda}:= \{ E_{f}| f \in B  \}.
\end{equation}
Then \eqref{spd of T} implies that the dual basis to $\mathcal{B}_{\lambda}$ under $\langle \cdot , \cdot \rangle_{X} $ is 

\begin{equation}
\mathcal{B}_{\lambda}^* = \{E_{M^{-1}(f^*)}| f^* \in B^* \}. 
\end{equation}

For $E \in \mathcal{S}(X)$, we may define the standard action of $W$ on $\mathcal{S}(X)$ as 
\begin{equation}
    w(E) := E(w^{-1}\cdot). 
\end{equation}

If $\mathcal{B}_{\lambda}$ is transitive under $W$, then the inner spectral coefficients of each basis element may be obtained by computing them for just one basis element and then applying $W$ to obtain all the other ones. If $B$ and $B^*$ are both transitive under $W$, then so are $\mathcal{B}_{\lambda}$ and $\mathcal{B}_{\lambda}^*$. The simplest example of such a basis $B$ is given by the standard basis 

\begin{equation}
    B := \{\delta_w | w \in W \}.
\end{equation}
We refer to the Weyl spline basis associated to $B$ as the \textit{standard Weyl spline basis}. $B$ is self dual under $\langle \cdot, \cdot \rangle_{\ell^2(W)}$, so the dual basis may be understood in terms of the action of $M^{-1}$ on $B$.\\

\textit{Example: A biorthogonal system of non-cardinal B-splines} \\

Let $\textbf{t} := \{t_0,t_1, \dots , t_{m}\}$ be a finite set of real numbers where $t_i < t_{i+1}$ for $0 \leq i\leq m-1$ and $m \geq 2$. We may recall that the B-spline $B_{\textbf{t}}$ is, up to scaling, the unique function in $C^{m-2}(\R)$ which is supported on $[t_0, t_m]$ and whose restriction to $[t_i,t_{i+1}]$ is a polynomial of degree at most $m-1$ for all $0 \leq i\leq m-1$.\\

Let $k$ be a positive integer, $R = \{-1,1\}$ be a one-dimensional root system, and  

\begin{equation}
X = [1, \dots, 1] \in \R^{2k}. 
\end{equation}

The Weyl group in this case consists of the identity map and the map $x \mapsto -x$. Let $B_{2k} := B_X$ which is the centered cardinal B-spline of order $2k$ and let $T_{2k} := T_X$ so that 
\begin{equation}
    T_{2k}(x) = \begin{cases}
        \frac{x^{2k-1}}{(2k-1)!} \ &  \text{for} \ x \geq 0, \\
        0 \ & \text{otherwise}.
\end{cases}
\end{equation}

For any positive integer $N$, the Weyl splines are, generically, supported on the closed interval $[-(N+k),(N+k)]$. Let $E_{N,k} \in \cS_X(N)$ be the Weyl spline satisfying 
\begin{equation}
    c_{E_{N,k}}(N) = 1 \ \text{and} \  c_{E_{N,k}}(-N) = 0
\end{equation}
and let $E_{-N,k}$ be the reflection of $E_{N,k}$. This is exactly the standard Weyl spline basis. It follows from this and Theorem \ref{GeneralizedWeylCharacter} that the $T_{2k}$-representation of $E_{N,k}$ is supported on  

\begin{equation}
   \textbf{t}_{N}(k) = \{\pm N, \pm (N+1), \dots , (N +k ) \}.
\end{equation}

The definition of the B-spline implies $E_{N,k} = B_{\textbf{t}_N(k)}$ and $E_{-N,k} = B_{-\textbf{t}_N(k)}$ (both equalities up to scaling). So the standard Weyl basis for $\mathcal{S}_X(N)$ consists of B-splines. The dual basis also consists of B-splines. The dual $E_{N,k}^*$ to $E_{N,k}$ under $\langle \cdot , \cdot \rangle_X$ satisfies that $D^{2k}E_{N,k}^*$ is supported on the $2k+1$ points   

\begin{equation}
   \textbf{t}^*_{N}(k) =  \{N, \pm (N+1), \dots \pm(N +k) \} 
\end{equation}

So, $E_{N,k}^* = B_{\textbf{t}^*_{N}(k)}$ and  $E_{-N,k}^* = B_{-\textbf{t}^*_{N}(k)}$ (both equalities up to scaling). With these Weyl splines being B-splines, they can be computed via classical B-spline theory (see, e.g., \cite{Splines}). 

\begin{figure}[h!]
\caption{Below are plots of the standard Weyl spline basis for $\cS_X(N)$ in the piece-wise linear setting. Equivalently, they are the linear B-spline with knot sequence $\{\pm N, (N+1) \}$ and its image under the reflection $x \mapsto -x$. }

\begin{subfigure}[h!]{0.4\linewidth}
\hspace{2.05cm}
\includegraphics[width = 14.1cm ,height = 3.3cm]{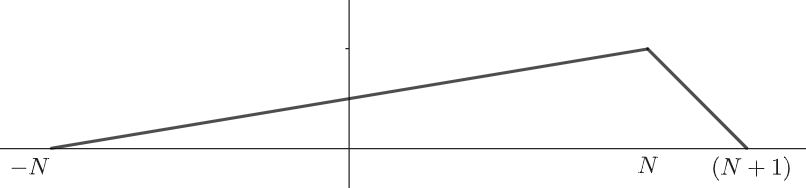}
\end{subfigure}

\begin{subfigure}[h!]{0.4\linewidth}
\vspace{1.5mm}
\includegraphics[width = 14.1cm ,height = 3.3cm]{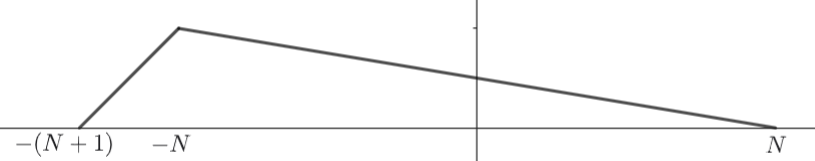}
\end{subfigure}
\end{figure}

\begin{figure}[h!]
\caption{Below are plots of the dual basis of the standard Weyl spline basis for $\cS_X(N)$ under $\langle \cdot , \cdot \rangle_X$ in the piece-wise linear setting which is plotted above. Equivalently, this is the linear B-spline with knot sequence $\{ N, \pm(N+1) \}$ and its image under the reflection $x \mapsto -x$. }

\begin{subfigure}[h!]{0.4\linewidth}
\includegraphics[width = 16.1cm ,height = 3.6cm]{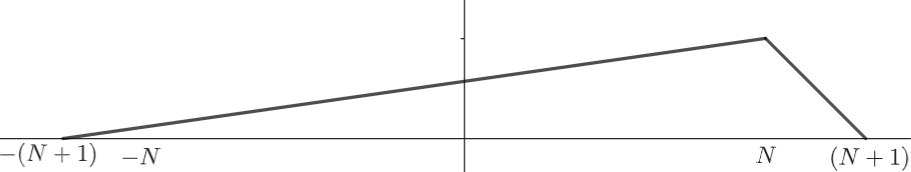} 
\end{subfigure}

\begin{subfigure}[h!]{0.4\linewidth}
\vspace{1.5mm}

\includegraphics[width = 16.1cm ,height = 3.7cm]{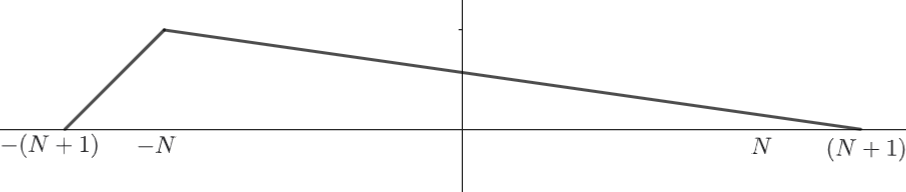}
\end{subfigure} 
\end{figure}

\newpage
We finish the section by proving the existence of one more basis for $\mathcal{S}_X(\lambda)$ with computationally desirable properties.  
\begin{corollary}

There exists an orthogonal basis for $\mathcal{S}_X(\lambda)$ which is transitive under $W$.
\end{corollary}

\begin{proof}
Let $E_f,E_g\in \calS_X(\lambda)$. Then 
\begin{equation}
\langle f, Mg \rangle_{\ell^2(W)} =\langle E_f, E_g \rangle_X  = \langle E_g, E_f \rangle_X = \langle g, Mf \rangle_{\ell^2(W)}, 
\end{equation}
hence $M$ is self-adjoint. It is positive definite (with respect to $\langle \cdot , \cdot \rangle_{\ell^2(W)}$) since $\langle \cdot, \cdot\rangle_X $ is an inner product. Thus, $M$ has a positive definite square root $M^{1/2}$ and since $M$ is invertible, so is $M^{1/2}$. $M$ is determined by the action of $D_X$ on $\mathcal{S}_X(\lambda)$ and since $\cS_X(\lambda)$ is $W$-invariant and $W$ commutes with $D_X$ on $\mathcal{S}(X)$, $M$ commutes with the action of $W$ on $\ell^2(W)$. It is well known that $M^{-1/2}$ is a polynomial in $M$, so $M^{-1/2}$ commutes with $W$. Therefore,
\begin{equation}
\overline{B} := \{M^{-1/2}\delta_w | w \in W\} 
\end{equation}
implies the existence of our desired orthogonal basis.
\end{proof}

\section{Appendix}

\subsection{The Weyl Character Formula and Box Splines Associated to Root Systems}
In this section, we briefly review root systems and explain how the classical Weyl character formula
yields the box-spline version of the Weyl character formula stated in Theorem \ref{splineversionofweyl}. For any proofs of the theorems below, we refer to \cite{Hall}. \\

A finite set of non-zero vectors $R \subset \R^d$ is a \textit{root system} if 
\begin{itemize}
    \item $R$ spans $\R^d$
    \item For each $\alpha \in R$, the only scalar multiples of $\alpha$ that belong to $R$ are $\pm\alpha$.
    \item For each $\alpha \in R$, $R$ is closed under the reflection about the hyperplane perpendicular to $\alpha$ i.e, for all $\beta \in R$
    \begin{equation}
        s_{\alpha}(\beta) = \beta - 2\frac{\alpha \cdot \beta}{||\alpha||^2}\alpha \in R.
    \end{equation}
    \item $2\displaystyle{\frac{(\alpha \cdot \beta)}{||\alpha||^2}}$ is an integer for all $\alpha , \beta \in R$. \\
\end{itemize}

$R$ is \textit{irreducible} if it cannot be partitioned into two non-empty mutually orthogonal subsets. \\ 

The Weyl group $W$ of $R$ is the group generated by reflections $\{s_{\alpha}| \alpha \in R \}$. The positive roots, $R^+$, is a subset of $R$ which satisfies the following properties 

\begin{enumerate}
    \item For each $\alpha \in R$, exactly one of $\alpha$ and $-\alpha$ lies in $R^+$. 
    \item For each $\alpha , \beta \in R^+$, if $\alpha + \beta \in R$, then $\alpha + \beta \in R^+$.
\end{enumerate}

The \textit{dominance ordering} $\leq$ is the following partial ordering on $\R^d$: $x \leq y$  if and only if 
 $$ y-x \in \text{cone}(R^+).$$
 
This partial ordering helps us distinguish between the distinct orbits of $\R^d$ under $W$. To see this, we let $H_{\alpha}$ be the hyperplane orthogonal to $\alpha$ for each $\alpha \in R$. The complement of 

\begin{equation}
    \bigcup_{\alpha \in R^+} H_{\alpha}
\end{equation}
is disconnected and we refer to each connected component as a \textit{Weyl chamber}.

\begin{theorem}
    The Weyl group acts freely and transitively on the Weyl chambers. Moreover, suppose $C$ is the (open) fundamental Weyl chamber i.e

\begin{equation}
    C :=\{v \in \R^d |   v \cdot \alpha > 0 \ \text{for all} \  \alpha \in R^+ \}
\end{equation}    
    and let $\overline{C}$ be the topological closure of $C$. For any $x \in \R^d$, there is exactly one element $x^*$ in the orbit $W(x)$ lying in $\overline{C}$.
\end{theorem}

Using the notation of this theorem, for any subset $S \subset \R^d$, we define $S^*$ to be
\begin{equation}
    S^* := \{x^*| x \in S \}.
\end{equation}

For any $x \in \R^d$, the \textit{Weyl polytope} associated with $x$ is 

\begin{equation}
    C_x := \text{conv}(W(x)).
\end{equation}
The following theorem shows the dominance ordering on $\R^d$ is related to whether certain Weyl polytopes contain each other.

\begin{theorem}
   For any $x, y \in C$, $x \leq y$ if and only if $x \in C_{y}$.
\end{theorem}

For a choice of positive roots $R^+$, its corresponding set of simple roots $\Pi$ is defined to be the subset of $R^+$ such that no element of $\Pi$ is a sum of two positive roots i.e
\begin{equation}
\Pi := \{\alpha \in R^+| \ \alpha \neq \beta + \gamma \ \text{for any $\beta,\gamma \in R^+$}  \}.
\end{equation}

The simple roots form a basis for $\R^d$ and satisfy that every $\alpha \in R^+$ is a linear combination of the simple roots with all coefficients being non-negative integers.  For any $x \in C$, the vertices of $C_{x}$ adjacent to $x$ are $\{s_{\alpha}(x)| \alpha \in \Pi \}$.  \\

We may now work towards Theorem \ref{splineversionofweyl}. For a root system $R$, we define the root lattice $\Lambda$ to be the lattice generated by the roots i.e 

\begin{equation}
    \Lambda := \bigg\{\sum_{\alpha \in R} c_{\alpha} \alpha \bigg| c_{\alpha} \in \Z \bigg\}.
\end{equation}

and we define the weight lattice  
\begin{equation}
   \Gamma := \bigg\{ \gamma \in \R^d \bigg| \frac{2(\gamma \cdot \alpha)} {||\alpha||^2} \in \Z \ \text{for all $\alpha \in R$} \bigg\}.
\end{equation} 

It is well known that the weight lattice $\Gamma$ contains the root lattice $\Lambda$ as a subgroup of finite index and is invariant under $W$. We refer to the elements of $\Gamma^*$ as the \textit{dominant weights}. \\ 

Recall the Weyl vector associated to a root system $R$ is $\overline{\rho} := \frac{1}{2}\rho_{R^+}. $

\begin{theorem}
    $\overline{\rho}$ lies in the open fundamental Weyl chamber.
\end{theorem}

We define the algebra of exponential polynomials on $\Gamma$ with real coefficients as 
\begin{equation}
     E(\Gamma) := \text{span}_{\R}\{e^{\gamma} | \gamma \in \Gamma \} 
\end{equation}


where $e^{\mu + \lambda} = e^{\mu}e^{\lambda}$ for any $\mu, \lambda \in \Gamma$. We have a natural action of $W$ on $E(\Gamma)$ which satisfies $w(e^{\mu}) = e^{w(\mu)}$. We define 

\begin{equation}
    E(\Gamma)_W = \{p \in E(\Gamma)| w(p) = p \ \text{for all } w \in W\}.
\end{equation}

Now, let $S_{\Pi} = \{s_{\alpha}| \alpha \in \Pi \}$ and define the length function $\ell$ on $W$ to be
 
\begin{equation}
    \ell(w) := \text{smallest number of elements in $S_{\Pi}$ needed to write $w$ as a product.}
\end{equation} 
The following theorem may be viewed as the \textit{Weyl Character Formula} stripped of its representation theoretic context.


\begin{theorem}\label{Weyl_char_form}
For each $\lambda \in \Gamma^*$, there exists a unique $\chi_{\lambda} \in E(\Gamma)_W$ satisfying
\begin{equation}\label{WeylCharFormula}
     \bigg(\prod_{\alpha \in R^+} (e^{\alpha/2} - e^{-\alpha/2}) \bigg) \chi_{\lambda} = \sum_{w \in W} (-1)^{\ell(w)}e^{w(\lambda + \overline{\rho})}.
\end{equation}
Moreover, $e^{\mu}$ has non-zero coefficient in $\chi_{\lambda}$ if and only if $\mu \in \Lambda + \lambda$ and $\mu^* \leq \lambda$.
\end{theorem}

It is worth noting an interpretation of this theorem which says that $\chi_{\lambda}$, which, generically speaking, is a large linear combination of exponentials, when multiplied by a particular collection of even more exponentials, results in a linear combination of exactly $|W|$ exponentials.\\

To connect this theory to box splines, we fix a root system $R$. Let $L$ be the linear map from $E(\Gamma)$ to $S(\Gamma)$ which sends 
\begin{equation}
    L: e^\mu \mapsto (-1)^{|R^+|}\delta(\cdot - \mu).
\end{equation}
for all $\mu \in \Gamma$. If $s_{\lambda}$ is the image of $\chi_{\lambda}$ under this map, the left-hand side of \eqref{WeylCharFormula} is mapped, under $L$, to the sequence corresponding to distribution
\begin{equation}
D_{R^+} (B_{R^+} \ast' s_{\lambda}).
\end{equation}

After convolving sequence with $T_{R^+}$, we obtain

\begin{equation}
B_{R^+} \ast' s_{\lambda} = \sum_{w \in W} (-1)^{\ell(w)}T_{R^+}(\cdot - w(\lambda + \overline{\rho})).
\end{equation}
and get Theorem \ref{splineversionofweyl}. Since the lattice generated by $R^+$ is $\Lambda$, the root lattice (not the weight lattice), we only make use of Theorem \ref{Weyl_char_form} in the case where $\lambda \in \Lambda^*$.

\subsection{Building $\cS_X(\lambda)$ from Opdam-Heckman Theory}
Here we give a very brief summary of the theory of orthogonal polynomials associated to root systems which we used throughout this paper. In particular, we show the construction of the \textit{non-symmetric Heckman-Opdam polynomials} (also referred to as \textit{non-symmetric Jacobi polynomials}). For proofs of these theorems, we refer to \cite{O}.\\

Define the coroot lattice $\overline{\Lambda}$ of $R$ to be the integer linear combinations of the \textit{coroots} i.e
\begin{equation}
\overline{\Lambda} := \text{span}_{\Z}\bigg\{\frac{2\alpha}{||\alpha||^2} \bigg| \alpha \in R \bigg\}
\end{equation}
Define the torus 
\begin{equation}
T:= \R^d/\overline{\Lambda}.    
\end{equation}

The map $e^{\gamma}: \xi \to e^{2\pi i \langle \gamma ,\xi \rangle }$ is an algebraic character of $T$ and $E(\Gamma)$ is the span of all algebraic characters of $T$. We let $dt$ be the Haar measure on $T$ that is normalized so that $\int_T dt= 1 $. Let $k:R \to \R_{\geq 0}$ be a $W$-invariant function, referred to as a \textit{multiplicity function}, and denote its values as $k_{\alpha}$. We may define the weight function:

\begin{equation}
    \delta_k := \prod_{\alpha \in R} (1-e^{\alpha})^{k_{\alpha}} 
\end{equation}

Consider the Hermitian inner product $(\cdot, \cdot )_k$ defined on $E(\Gamma)$ as

\begin{equation}
    (f,g)_k :=  \int_T f(t) \overline{g(t)}  \delta_k(t) \ dt.
\end{equation}

Let $\leq_B$ denote the \textit{Bruhat ordering} on $W$ defined as $w \leq_B v$ if some sub-string (not necessarily consecutive) of some reduced word for $v$ is a reduced word for $w$. Let $W_{\lambda}$ be the stabilizer of $\lambda$ under $W$ and define the set
\begin{equation}
    W^{\lambda} := \{w \in W| \ell(w w') \geq \ell(w) \ \text{for all} \ w' \in W_{\lambda} \}.
\end{equation}

It is well known that $W^{\lambda}$ is a set of coset representatives of $W_{\lambda}$ in $W$. Finally, let $w^{\lambda} \in W^{\lambda^*}$ be the unique element such that $\lambda = w^{\lambda}\lambda^*$. The following is a partial ordering on $\Gamma$: 

\begin{equation}
\gamma \leq_W \mu \Longleftrightarrow
\begin{cases}
    (1) & \gamma^* < \mu^*  \ \text{if}  \ \gamma^* \neq \mu^*  \\
    (2) & w^{\gamma} \leq_B w^{\mu} \ \text{if $\gamma^* = \mu^*$}
\end{cases}
\end{equation}

where $\gamma <_W \mu$ if and only if $\gamma \leq_W \mu$ and $\gamma \neq \mu$. We may now define the non-symmetric Heckman-Opdam polynomials. We define $E_k(\lambda)$ as the unique element in $E(\Gamma)$ satisfying the following two conditions:

\begin{enumerate}
    \item $E_k(\lambda) = \displaystyle{e^{\lambda} + \sum_{\mu <_W \lambda} c_{\mu\lambda}e^{\mu} } $
    \item $(E_k(\lambda), e^{\mu} )_{k} = 0 $ for all $\mu  <_W \lambda$.
\end{enumerate}


\begin{theorem}\label{Opdams Main Theorem}
    $\{E_k(\lambda) \}_{\lambda \in \Gamma}$ is an orthogonal basis for $E(\Gamma)$ under $(\cdot, \cdot)_k$.
\end{theorem}


We take several steps in order to fit these results to a spline-theoretic context. Recall that $\Gamma^* = \overline{C} \cap \Gamma$ where $C$ is the fundamental Weyl chamber. From here on, we assume $\lambda \in \Gamma^*$ and define the space $S_k(\lambda)$ to be the subspace of $E(\Gamma)$ such that all $E \in S_k(\lambda)$ satisfy

\begin{enumerate}
    \item $E = \sum_{\mu^* \leq \lambda} c_{\mu} e^{\mu} $
    \item $(E, e^{\mu} )_k = 0 $ for all $\mu^* < \lambda$.
\end{enumerate}

By Theorem \ref{Opdams Main Theorem} and the definition of the non-symmetric Heckman-Opdam polynomials, we may conclude that  
 
\begin{equation}
  S_k(\lambda) = \text{span}\{E_k(\mu) \}_{\mu \in W(\lambda)} 
\end{equation}
and may then further conclude that the following decomposition is orthogonal:

\begin{equation}
    E(\Gamma) = \bigoplus_{\lambda \in \Gamma^*} S_k(\lambda).
\end{equation}

Now observe that since $k$ is a $W$-invariant function, we have that 
\begin{equation}\label{W_invariant_function}
    \sum_{\alpha \in R} k_{\alpha}\alpha = 0
\end{equation}

From here on, we suppose that the image of the multiplicity function $k$ lies in the set of positive integers. We may then rewrite $\delta_k$ using \eqref{W_invariant_function} as follows:   

\begin{align}
    \delta_k = \prod_{\alpha \in R} (e^{\alpha/2} - e^{-\alpha/2} )^{k_{\alpha}}
\end{align}

We may now associate a multi-set $X_0$ to $k$ as defined in Section $1.3$. We have that
\begin{align}
    \prod_{\alpha \in R} (e^{\alpha/2} - e^{-\alpha/2} )^{k_{\alpha}} = (-1)^{|X_0|}  \prod_{\alpha \in R^+} (e^{\alpha/2} - e^{-\alpha/2} )^{2k_{\alpha}}.
\end{align}

The map from $E(\Gamma)$ to $S(\Gamma)$ which sends  

$$e^{\mu} \to \delta(\cdot - \mu),$$ 
for all $\mu \in \Gamma$, sends $\delta_k$ to a non-zero scalar multiple of $\Delta_X$ where $X = [X_0|X_0]$. Now, we restrict ourselves to the case where $\mu \in \Lambda$ and consider the linear map which sends

\begin{equation}
    e^{\mu} \mapsto B_X(\cdot-\mu)
\end{equation}

for any $\mu \in \Lambda$. This map, up to a non-zero scalar multiple, is an isometry between the sub-algebra of $E(\Gamma)$  
\begin{equation}
    E(\Lambda) := \text{span}_{\R}\{ e^{\mu} | \mu \in \Lambda \}
\end{equation}
and $\mathcal{S}(X)$ where $E(\Lambda)$ is endowed with $(\cdot , \cdot)_k$ and $\mathcal{S}(X)$ with $\langle \cdot, \cdot \rangle_X$. For any $\lambda \in \Lambda^*$, the image of $S_k(\lambda)$, under this isometry, is precisely $\cS_X(\lambda)$. \\

\subsection{Local Discrete Quasi-Interpolants Can Be Made Optimal}
 In this section, we write out the proof that a sequence $a \in S(\Lambda)$ can be chosen to make $Q_a$ reproduce polynomials. Carl de boor gives several ways of doing this in \cite{deBoor1990}. Thus, we assume the notation throughout Section $3$ and give one argument here for completeness.

\begin{proof}[Proof of Theorem~\ref{Section 3.2 Main Result}]
Recall that a polynomial $p \in \Pi_d$ is determined by its values on $\Lambda$. Moreover, the semi-discrete convolution of $B_X$ with $p|_{\Lambda}$, where $p$ is of order $\kappa(X)$, is also a polynomial of order $\kappa(X)$ (see \cite{BoxSplines}). Thus, the desired polynomial reproduction property of $Q_a$ is equivalent to the statement 
\begin{align}
    Q_a(p)|_{\Lambda} = p|_{\Lambda}.
\end{align}

With this in mind, we observe now that any finitely supported sequence on $\Lambda$ is a linear combination of difference operators. This folows from the fact that for any $\mu \in \Lambda$, we have 

\begin{equation}
    \delta(\cdot - \mu) = \delta + \frac{1}{2}\bigg(\Delta^2_{\mu} - \Delta_{2\mu}\bigg). 
\end{equation}

Moreover, $\Lambda$ is generated by integer linear combinations of the simple roots $\Pi$. Thus, $S(\Lambda)$ is finitely generated as an algebra over the real numbers with generating set

\begin{equation}
    \{\Delta_{\alpha}^2 , \Delta_{2\alpha} | \alpha \in \Pi  \}.
\end{equation}

If we label the difference operators in this generating set as $\Delta_1 , \dots , \Delta_{2d}$, then there exists a polynomial $q \in \Pi_{2d}(\R)$ such that
\begin{equation}
    B_X|_{\Lambda} = q(\Delta_{1}, \dots , \Delta_{2d}).
\end{equation}

As a result, $q(x)$ has a non-zero constant term and we may write 

\begin{equation}
    B_X|_{\Lambda} = \frac{\delta}{\text{det}(\Lambda)} + g(\Delta_1, \dots ,\Delta_{2d})
\end{equation}

where $g \in \Pi_{2d}(\R)$ has constant term $0$. It is not difficult to construct a polynomial $a(x) \in \Pi_{2d}(\R)$ such that 
\begin{equation}
    q(x)a(x) = 1 + r(x) 
\end{equation}
where the monomial of lowest degree in $r(x)$ is at least $\kappa(X)+1$. Now recall that for any multi-set $S \subset \R^d$ containing $N$ vectors, any polynomial of order $N$ in $\Pi_d(\R)$ lies in the kernel of $\Delta_S$. Thus, if $p \in \Pi_d(\R)$ is of degree $\kappa(X)$, we must have 
\begin{align}
    r(\Delta_1, \dots \Delta_{2d})p = 0
\end{align}
since every term in $r(\Delta_1, \dots \Delta_{2d})$ is a difference operator of order at least $\kappa(X)$. If we let
\begin{equation}
    a := a(\Delta_1, \dots ,\Delta_{2d}),
\end{equation}

then we have
\begin{align}
    (B_X|_{\Lambda} \ast a) \ast p|_{\Lambda} & = (\delta + r(\Delta_1, \dots \Delta_N)) \ast p|_{\Lambda} \\
    & = p|_{\Lambda}
\end{align}
\end{proof}

\centering\noindent\textbf{Acknowledgments.}\\
The author was supported in part by NSF DMS- 2037851 \\

AI Statement: The author used artificial-intelligence tools solely for language editing, limited assistance with local mathematical reasoning, and minor image editing. All mathematical content was independently written and checked by the author, who takes full responsibility for it.
\bibliographystyle{unsrtnat}
\bibliography{main}

\bigskip

\noindent
\textsc{University of Wisconsin–Madison, Department of Mathematics, 480 Lincoln Drive, Madison, WI 53706 }

\noindent
\textit{Email address:}
\texttt{djaime@wisc.edu}
\end{document}